\hsize = 31pc
\vsize = 45pc
\input amssym.def
\input amssym.tex
 \font\newrm =cmr10 at 24pt
\def\bul{\raise .9pt\hbox{\newrm .\kern-.105em } }

 \def\fr{\frak}

\baselineskip=13pt
 
 \def\h{\hbox{ }}
 
 \def\p{{\fr p}}
 \def\u{{\fr u}}

 \def\m{{\fr m}}
 \def\n{{\fr n}}
 \def\a{{\fr a}}

 \def\ss{{\fr s}}
 \def\k{{\fr k}}

 \def\hh{{\fr h}}
 \def\tt{{\fr t}}

 \def\g{{\fr g}}
 \def\v{{\fr v}}
 
 \def\q{{\fr q}}
 
 \def\x{{\chi}}

 \def\<{\le}
 \def\>{\ge}

 \def\s{{\h\subset\h}}
 
 \def\vs{\vskip }

 \def\mapright#1
  {\smash{\mathop
  {\longrightarrow}
  \limits^{#1}}}

 \def\kk#1{{\kern .4 em} #1}
 \def\vs{\vskip 1pc}

\font\twelverm=cmbx12 at 14pt
\hsize = 31pc
\vsize = 45pc
\overfullrule = 0pt

\def\vs{\vskip 1pc}
\font\ninerm=cmr10 at 9pt

\font\smallbf=cmbx10 at 9pt 

\font\twelverm=cmbx12 at 14pt
\font\authorfont=cmbx10 at 12pt
%\font\ninerm=cmr9
\centerline{\twelverm Minimal coadjoint orbits and symplectic induction} 
\vskip 1.5pc
\baselineskip=11pt
\vskip8pt
\centerline{\authorfont BERTRAM KOSTANT\footnote*{\ninerm
Research supported in part by NSF contract DMS-0209473 and the KG\&G Foundation.}}
\vskip 2pc
\line{\hfil{\it To Alan, with admiration, on the occasion of his sixtieth birthday}}

\vskip 1.5pc
{\font\ninerm=cmr10 at 9pt \baselineskip=14pt 
\noindent{\smallbf ABSTRACT.} \ninerm Let $(X,\omega)$ be an integral symplectic
manifold and let $(L,\nabla)$ be a quantum line bundle, with connection, over $X$ having
$\omega$ as curvature.  With this data
one can define an induced symplectic manifold $(\widetilde {X},\omega_{\widetilde
{X}})$ where $dim\,\widetilde {X} = 2 + dim\,X$. It is then shown that prequantization
on $X$ becomes classical Poisson bracket on $\widetilde {X}$. We consider the
possibility that if $X$ is the coadjoint orbit of a Lie group $K$ then $\widetilde
{X}$ is the coadjoint orbit of some larger Lie group $G$. We show that this is the case
if $G$ is a non-compact simple Lie group with a finite center and $K$ is the maximal
compact subgroup of $G$. The coadjoint orbit $X$ arises (Borel-Weil) from the action of
$K$ on
$\p$ where $\g= \k +\p$ is a Cartan decomposition. Using the Kostant-Sekiguchi
correspondence and a diffeomorphism result of M. Vergne we establish a symplectic
isomorphism
$(\widetilde {X},\omega_{\widetilde {X}})\cong (Z,\omega_Z)$ where $Z$ is a
non-zero minimal ``nilpotent" coadjoint orbit of $G$. This is applied to show that the
split forms of the 5 exceptional Lie groups arise symplectically from the symplectic
induction of coadjoint orbits of certain classical groups.} 

\vskip 1pc
\baselineskip 15pt
\centerline{\bf 0. Introduction}\vskip 1pc

\rm
 0.1. Let $(X,\omega)$ be a connected symplectic manifold and let $Ham(X)$ be
the Lie algebra of all smooth Hamiltonian vector fields on $X$. The
space $C^{\infty}(X)$ of all smooth $\Bbb C$-valued functions
on
$X$ is a Lie algebra under Poisson bracket. To any  
$\varphi\in C^{\infty}(X)$ there corresponds $\xi_{\varphi}\in
Ham(X)$ and $\varphi\mapsto \xi_{\varphi}$ realizes $C^{\infty}(X)$
as a Lie algebra central extension $$0\longrightarrow \Bbb
C\longrightarrow  C^{\infty}(X) \longrightarrow Ham(X)\longrightarrow 0\eqno (0.1)$$
of $Ham(X)$ by the constant functions.

Prequantization, when it exists, is a specific representation of the Poisson Lie algebra
$C^{\infty}(X)$ which does not descend to $Ham(X)$ (Heisenberg-like --- it is
non-trivial on the constant functions). A necessary and sufficient condition for
prequantization is that the deRham class $[\omega]\in H^2(X,\Bbb R)$ lie in the image of
the natural map $H^2(X,\Bbb Z)\to H^2(X,\Bbb R)$. In such a case we will say
$(X\,\omega)$ (or just $X$ if $\omega$ is understood) is integral. If $(X\,\omega)$ is
integral there exists  a complex line bundle
$L$ (the quantum line bundle) with connection
$\nabla$ over $X$ such that
$$\omega = curv\,(L,\nabla)\eqno (0.2)$$ Using the connection one defines the
covariant derivative $\nabla_{\xi}\,s$ of any smooth section $s$ of $L$ by any vector
field (v.f.) on $X$. Furthermore the connection can (and will) be chosen so that there
exists a Hilbert space structure in each fiber of $L$ which is invariant under
parallelism. By considering only the unit circles in each fiber of $L$ one obtains a
principal $U(1)$-bundle $$ \matrix{ U(1) & \kern-1em\longrightarrow &\kern-1em
L^1 &\cr
										& &\kern-1.4em \big\downarrow &\kern-2.6em\tau\cr
& && \kern-3.6em X\cr
\noalign{\vskip3pt}\cr}\eqno (0.3)$$ The connection defines a real
$U(1)$-invariant 1-form $\alpha$ on $L^1$ which on each fiber corresponds to
$d\theta/2\pi$ on $U(1)$ and one has $$d\alpha = \tau^*(\omega)\eqno (0.4)$$

Let $S$ be the linear space of all smooth sections of $L$. Then prequantization is
the Lie algebra representation $\pi$ of $C^{\infty}(X)$ on $S$ given by
$$\pi(\varphi)\,s =
(\nabla_{\xi_{\varphi}} + 2\,\pi\,i\,\varphi)\,s\eqno (0.5)$$ Let $\zeta$ be the
vertical vector field on $L^1$ (generating the $U(1)$-action) such that
$\langle\alpha,\zeta\rangle = -1$. One has a linear isomorphism $$S\to \widetilde
{S}\s C^{\infty}(L^1),\qquad s\mapsto \widetilde {s}$$ where $\widetilde {S} = \{f\in
C^{\infty}(L^1)\mid
\zeta\,f = 2\,\pi\,i\,f\}$ and a (associative) algebra isomorphism $$C^{\infty}(X)\to
\widetilde {C}\s C^{\infty}(L^1),\qquad\varphi\mapsto \widetilde {\varphi}$$ where
$\widetilde {C} = \{f\in
C^{\infty}(L^1)\mid
\zeta\,f = 0\}$. 

Now let $\widetilde {X} = L^1\times \Bbb R^{+}$ so that $$dim\,\widetilde {X} = 2 +
dim\,X\eqno (0.6) $$ Let $r\in C^{\infty}(R^+)$ be the natural coordinate function on
$\Bbb R^+$ so that if $t\in \Bbb R^+$ then $r(t) = t$. One defines a symplectic
form $\omega_{\widetilde {X}}$ on $\widetilde {X}$ by putting $\omega_{\widetilde {X}}
= d\,(r\,\alpha)$ so that $$\omega_{\widetilde {X}} = d\,r\wedge \alpha +
r\,\widetilde {\omega}\eqno (0.7)$$ where we have put $\widetilde {\omega}= d\alpha$.
One notes that $\xi_{r} = \zeta$ and hence $(X,\omega)$ arises from $(\widetilde
{X},\omega_{\widetilde {X}})$ by (Marsden-Weinstein) symplectic reduction on the
hypersurface $r=1$. Reversing the direction we refer to the construction of
$(\widetilde {X},\omega_{\widetilde {X}})$ from $(X,\omega)$ as symplectic induction.
Among the statements in the following theorem is the result that prequantization in
$X$ is Poisson bracket in $\widetilde {X}$. In a word, quantization, at least at the
prequantized level, is classical mechanics two dimensions higher. \vs {\bf Theorem 0.1.}
{\it The map $$C^{\infty}(X)\to C^{\infty}(\widetilde X),\qquad \varphi\mapsto r
\widetilde
\varphi\eqno (0.8)$$ is a monomorphism of Poisson Lie algebras. Moreover for any $s\in
S$ and
$\varphi\in C^{\infty}(X)$ one has
$$\widetilde {\pi(\varphi)(s)} = [r\,\widetilde \varphi, \widetilde s]\eqno (0.9)$$}
See Theorem 1.6.

0.2. If $(X,\omega)$ is the coadjoint orbit of some Lie group $K$, where $\omega$ is
the KKS-symplectic form, consider the possibility that 
$(\widetilde {X},\omega_{\widetilde {X}})$ is the coadjoint orbit of some larger Lie
group $G$. It will be the main theorem of this paper to show that this indeed is
true in an important specialized case and that in this case $(\widetilde
{X},\omega_{\widetilde {X}})$ is a minimal coadjoint orbit of $G$.

Assume that $K$ is a compact connected Lie group. If $V$ is a finite dimensional,
complex irreducible module for K (and hence for its complexification $K_{\Bbb C}$)
then by the Borel-Weil theorem there corresponds to $V$ an integral coadjoint orbit
$X(V)\s\k^*$ where $k=Lie\,K$ and $\k^*$ is the dual to $\k$. In fact $X =X(V)$ is an
isomorphism to the unique closed $K_{\Bbb C}$-orbit in the projective space $Proj\,V$
and $L$ is defined by considering the cone over this orbit in $V$. The orbit is the
projective image of the affine variety of extremal weight vectors in $V$.

Now assume that $G$ is a non-compact Lie group with finite center such that $\g =
Lie\,G$ is simple and that $K$ is a maximal compact subgroup of $G$. Then $G/K$ is a
non-compact symmetric space and one has a Cartan decomposition $\g =\k + \p$. The
complexification
$\p_{\Bbb C}$ is a $K$ (and hence $K_{\Bbb C}$) module via the adjoint representation. 
There are 2 cases to be considered. We will say that
$\g$ is of non-Hermitian type if $G/K$ is non-hermitian and $\g$ is of hermitian type if
$G/K$ is Hermitian. Let $I$ be an index parameterizing the irreducible $K$-submodules
$V^i,\,i\in I$ of $\p_{\Bbb C}$. In the non-Hermitian case $I$ has 1 element
(i.e. $\p_{\Bbb C}$ is irreducible ) and $I$ has 2 elements in the Hermitian case. In
any case let $X^i = X(V^i),\,\,i\in I$, so that $X^i$ is an integral coadjoint orbit of
$K$. Let $X=X^i,\,i\in I$. In the Hermitian case $X\neq -X$ and
$\{X,-X\}=
\{X^i\},\,i\in I$. In the non-Hermitian case $X = -X$. 

Now the affine variety of
extremal weight vectors $E^i$ in $V^i$ is a $K_{\Bbb C}$-orbit and, by the
Kostant-Sekiguchi theorem, corresponds to a ``nilpotent" $G$-coadjoint orbit $Z^i$ in the
dual space $\g^*$. If $Z = Z^i$ then $Z \neq -Z$ in the Hermitian case and $\{Z,-Z\} =
\{Z^i\},\,i\in I$. In the non-Hermitian case $Z = -Z$ and in any case if $Y$ is any
non-zero coadjoint $G$-orbit then $$dim\,Y \geq dim\, Z\eqno (0.10)$$ On the other
hand $C^{\infty}(Z)$ is a Lie algebra under Poisson bracket with respect to
the KKS symplectic form $\omega_Z$ on $Z$ and since $Z$ is a coadjoint orbit one has a
Lie algebra embedding $$\g\to C^{\infty}(Z)\eqno (0.11)$$ It follows from (0.10) and
the orbit covering theorem that $dim\,Z$ is the smallest possible dimension of a
symplectic manifold which has $\g$ as a subalgebra of functions under Poisson
bracket. A theorem of Mich\`ele Vergne asserts that Kostant-Sekiguchi corresponding
orbits in general are $K$-diffeomorphic. In the present case because of the minimality
(0.10) Vergne's diffeomorphism can be given very simply and it leads to a
$K$-diffeomorphism $$\widetilde {X}\to Z\eqno (0.12)$$ In particular $$dim\,Z = 2 +
dim\,X\eqno (0.13)$$ The following is our main theorem. In effect it says that
symplectically inducing the $K$-coadjoint orbit $X$ ``sees" the non-compact simple Lie
algebra $\g$. \vs {\bf Theorem 0.2.} {\it The map (0.12) now written
$$\beta:(\widetilde {X},\omega_{\widetilde {X}})\to (Z,\omega_Z) \eqno (0.14)$$ is a
symplectic diffeomorphism so that one has a (minimal) Lie algebra injection $$\g\to
C^{\infty}(\widetilde {X})\eqno (0.15)$$ Furthermore if $\mu$ is the moment map with
respect to the action of $K$ on $Z$ then $$\mu(Z) = \Bbb R^+\,X\eqno (0.16)$$}\vs  See
Theorems 3.10, 3.13 and 3.16. \vs {\bf Remark 0.3.} Note that (0.15) points to an
interesting difference between
$(X,\omega)$ and the induced symplectic manifold $(\widetilde {X},\omega_{\widetilde
{X}})$ in that one cannot have a non-trivial Lie algebra homomorphism $\g \to
C^{\infty}(X)$ since $dim\,X = dim\,\widetilde {X}\,-2$. Indeed the statement above
concerning
$Z$ now implies that $dim\,\widetilde {X}$ is
the smallest possible dimension of a symplectic manifold $W$ which admits an embedding
of
$\g$ as a Lie algebra of functions on $W$ under Poisson bracket.\vs 0.3. We say that $\g$
is 
$O_{min}$-split if $\g_{\Bbb C}$ is a simple complex Lie algebra and
$e\in O_{min}$ where $O_{\min}$ is the minimal nilpotent orbit in
$\g_{\Bbb C}$. The definition of $O_{min}$-split depends only on $\g$ and is, in 
particular, independent of the choice of $e$. The following result is Theorem 3.17. \vs
{\bf Theorem 0.4.}  {\it The simple Lie algebra $\g$ is $O_{min}$-split if and only if
$$dim\,Cent\,\n =1\eqno (0.17)$$} \vs {\bf Remark 0.5.} Another criterion for the
$O_{min}$-split condition was cited on the top of p. 19 in [B-K]. \vs Proposition 3.19
is stated here as \vs {\bf Proposition 0.6.} {\it
If $\g$ is $O_{min}$-split (e.g. if $\g$ is split) then $X$ is not only a
$K$-symmetric space but in fact $X$ is a Hermitian symmetric space.}\vs Assume that $\g$ is
$O_{min}$-split and $\g$ is of non-Hermitian type so that $\k_{\Bbb C}$ is
semisimple. Let $K'$ be a non-compact real form of $K_{\Bbb C}$ having $K_{\nu}$ as a
maximal compact subgroup so that
$K'/K_{\nu} = X'$ is the non-compact symmetric dual to the compact symmetric space
$X$. One can then show that not only is $X'$ a complex bounded domain but in fact
$X'$ is a tube domain. In particular by the Kantor-Koecher-Tits theory
$X'$ corresponds to a formally real Jordan algebra $J(X)$.

If $\g$ is a split form of any one of 5 exceptional simple Lie
algebras, $\g$ is non-Hermitian so that the statement above applies to $\g$. But a
minimal (dimensional) symplectic realization of $\g$ as functions on a symplectic
manifold is achieved when the manifold is the induced symplectic $\widetilde {X}$ of a
coadjoint orbit of a compact Lie group $K$. The group $K$ turns out to
be classical in all 5 cases so that the exceptional Lie algebras $\g$ emerge
symplectically from the symplectic induction of a classical coadjoint orbit. In \S
3.4 we present a table which contains the relevant information. The cases of
$E_6,\,E_7$ and
$E_8$ are taken from [B-K].\vs 0.4. We wish to thank Ranee Brylinski for many
conversations on related matters at an earlier time. A number of ideas in this paper
evolved when [B-K] was written. This is particularly true of Theorem 1.7 which therefore
should be considered collaborative.
\vs
 
\centerline{\bf 1. Symplectic Induction and the construction of $(\widetilde X,\omega_{\widetilde
X})$}\vskip 1.5pc
\rm 1.1. In this section we recall the theory of prequantization. See [K-1]. Let
$(X,\omega)$ be a connected symplectic manifold. For any
$\varphi\in C^{\infty}(X)$ one defines a Hamiltonian vector
field $\xi_{\varphi}$ on $X$ so that for any vector field
(v.f.)
$\eta$ one has $$\eta\,\varphi = \omega(\xi_{\varphi},\eta)$$
Poisson bracket in  $C^{\infty}(X)$ is defined by putting $[\varphi,\psi] =
\xi_{\varphi} \psi$ for any $\varphi,\psi\in C^{\infty}(X)$. If
$Ham(X) = \{\xi_{\varphi}\mid \varphi\in C^{\infty}(X)\}$ is the Lie
algebra of all Hamiltonian vector fields on $X$ then
$C^{\infty}(X)$, under Poisson bracket, is a Lie algebra central extension $$0\longrightarrow \Bbb
C\longrightarrow  C^{\infty}(X) \longrightarrow Ham(X)\longrightarrow 0\eqno (1.1)$$ of $Ham(X)$, by the
constant functions on $X$, via the map $\varphi\mapsto \xi_{\varphi}$. 

We now assume that the de Rham class $[\omega]\in H^2(X,\Bbb R)$ is integral (i.e., $[\omega]$ is in the
image of the natural map $H^2(X,\Bbb Z)\to H^2(X,\Bbb R)$). Then one knows that there exists a complex
line bundle
$L$ with connection
$\nabla$, over
$X$, such that
$$\omega = curv\,(L,\nabla) \eqno (1.2)$$ Using the connection one defines the covariant derivative
$\nabla_{\xi}\,s$ of any section $s$ of $L$ (either local or global, but always assumed to be infinitely
smooth) by any v.f.
$\xi$ on $X$. If $L^*$ is $L$ minus the zero
section then $L^*$ is a principal $\Bbb C^*$ bundle $$\matrix{\Bbb C^* & \kern-.9em\longrightarrow
&\kern-1em L^*
\cr
										& &\kern-1.5em \big\downarrow \cr
& & \kern-1.65em X\cr
\noalign{\vskip3pt}\cr}
\eqno (1.3)$$ In addition there exists a $\Bbb C^*$-invariant 1-form $\alpha$ on $L^*$ which, on any
fiber of (1.3), pulls back to ${1\over 2\,\pi\,i}\,{dz\over z}$ on $\Bbb C^*$ and is such that for any
local section $s$ of (1.3), and any v.f. $\xi$ on $X$, one has $$\nabla_{\xi}\,s = 2\,\pi\,i \langle
s^*(\alpha),\xi\rangle\,s\eqno (1.4)$$ on the domain of $s$ and on this domain $$d\,(s^*(\alpha)) =
\omega\eqno (1.5)$$\vskip .5pc The connection $\nabla$ may be (and will be) chosen so that there exists a
Hilbert space structure on each (1 dimensional) fiber of $L$ which is invariant under parallel
transport. See Proposition 2.1.1  in [K-1]. One then defines a principal $U(1) = \{e^{i\,\theta}\mid
\theta\in
\Bbb R\}$ bundle
$L^1$ over $X$ with bundle projection $\tau$, $$\matrix{ U(1) & \kern-1.2em\longrightarrow &\kern-1em
L^1 &\cr
										& &\kern-1.4em \big\downarrow &\kern-2.6em\tau\cr
& && \kern-3.6em X\cr
\noalign{\vskip3pt}\cr}\eqno (1.6)
$$ by taking the fibers of $L^1$ to be the unit circles in the corresponding fibers of $L$.
The restriction of $\alpha$ to $L^1$ is real and  on each fiber of $L^1$, pulls back, via (1.6), to the
1-form ${d\theta\over 2\,\pi}$ on $U(1)$. Henceforth we will identify $\alpha$ with this restriction and
(1.5) implies $$d\,\alpha =\tau^*(\omega)$$ Obviously $$dim\,L^1 = 1 + dim\,X\eqno
$$\vskip .5pc
1.2. Let $S$ be the space of all smooth global sections of $L$. Let $$\pi:C^{\infty}(X)\to End\,S$$ be
defined by putting, for any $\varphi\in C^{\infty}(X)$ and any $s\in S$, $$\pi(\varphi)(s) =
(\nabla_{\xi_{\varphi}} + 2\,\pi\,i\,\varphi)\,s\eqno (1.7)$$ The correspondence $\varphi\mapsto
\pi(\varphi)$ is called prequantization and it is a result (see Theorem 4.3.1 in [K-1]) that $\pi$ is a
Lie algebra representation of
$C^{\infty}(X)$ on $S$. (Theorem 4.3.1 in [K-1] is stated for real-valued functions on $X$. It
trivially extends to $C^{\infty}(X)$ by complex linearity.) Recalling (1.1) one notes that
$\pi$ does not descend to a representation of
$Ham(X)$.

Now with respect to the multiplication action of $U(1)$ on $\Bbb C$ we may 
write $L = L^1\times_{U(1)}\Bbb C$ so that $L$ is associated to the principal bundle $L^1$. If one
wishes $S$ may then be identified with the subspace $\widetilde {S}\s C^{\infty}(L^1)$ defined by
putting $$\widetilde {S} = \{f\in C^{\infty}(L^1)\mid f(q\,c) = c^{-1}\,f(q),\,\,\forall q\in L^1,\,c\in
U(1)\}\eqno (1.8)$$ The image of $s\in S$ under the linear isomorphism $$S\to \widetilde {S}\eqno
(1.9)$$ will be denoted by $\widetilde {s}$. 

A vector field on $L^1$ will be called vertical if it is tangent to the fibers
of
$L^1$. One notes that there there exists a unique (real) vertical field $\zeta$ on $L^1$ such that
$$\langle
\alpha,\zeta\rangle = -1\eqno (1.10)$$ With respect to an isomorphism of $U(1)$ with a fiber of $L^1$
note that the restriction of
$\zeta$ to that fiber corresponds to $-2\,\pi {d\over d\theta}$ on $U(1)$. Clearly the subspace
$\widetilde {S}$ may be also given by $$\widetilde {S} = \{f\in C^{\infty}(L^1)\mid \zeta f=
2\,\pi\,i\,f\}\eqno (1.11)$$ Now for any $\varphi\in C^{\infty}(X)$ let $\widetilde \varphi \in
C^{\infty}(L^1)$ be given by the pullback $\widetilde \varphi =\tau\circ \varphi$ and let $\widetilde
C = \{\widetilde \varphi\mid \varphi\in C^{\infty}(X)\}$. Clearly $$\widetilde
C = \{f\in C^{\infty}(L^1)\mid \zeta\,f = 0\}\eqno (1.12)$$ A vector field on $L^1$ will be called
horizontal if it is orthogonal to $\alpha$. If $\xi$ is any v.f. on $X$ it is clear that there exists a
unique horizontal vector field $\widetilde \xi$ on $L^1$ such that $\tau_*(\widetilde \xi) = \xi$. Now
for any $\varphi\in C^{\infty}(X)$ and vector field $\xi$ on $X$ let $\eta_{(\varphi,\xi)}$ be
the v.f. on $L^1$ defined by putting $$\eta_{(\varphi,\xi)} =  \widetilde \xi + \widetilde {\varphi}
\,\zeta
\eqno (1.13)$$ One readily establishes (see Proposition 2.9.1 in [K-1] for the case where $\Bbb C^*$
is used instead of $U(1)$)
\vs {\bf Proposition 1.1.} {\it The vector field
$\eta_{(\varphi,\xi)}$ on $L^1$ is $U(1)$-invariant and any $U(1)$-invariant v.f. on $L^1$ is uniquely
of this form.}\vs Now consider the question as to whether or not $$\theta(\eta_{(\varphi,\xi)})(\alpha)
=0\eqno (1.14)$$ where $\theta(\eta_{(\varphi,\xi)})$ is Lie differentiation by $\eta_{(\varphi,\xi)})$.
The following result interprets prequantization as the ordinary action
on $\widetilde S$ by those vector fields on $L^1$ which are (1) $U(1)$-invariant and (2) annihilate
$\alpha$ by Lie differentiation. For any $\varphi\in C^{\infty}(X)$ let $$\eta_{\varphi}=
\eta_{(\varphi,\xi_{\varphi})}\eqno (1.15)$$ Except for the use of $\Bbb C^*$ and $L^*$ instead of $U(1)$
and $L^1$ the following result is Theorem 4.2.1 in [K-1]. \vs {\bf Theorem 1.2.} {\it Let
$\varphi\in C^{\infty}(X)$. Then the vector field $\eta_{\varphi}$ on $L^1$ is (1) $U(1)$-invariant and
(2) annihilates $\alpha$ by Lie differentiation. Furthermore any v.f. on $L^1$ which satisfies (1) and
(2) is uniquely of this form. Moreover for any $s\in S$ one has $$\widetilde {\pi(\varphi)(s)} =
\eta_{\varphi}\widetilde {s}\eqno (1.16)$$} \vs 1.3. Let $\Bbb R^+$ denote the multiplicative group of
positive real numbers. Let $$\widetilde X = L^1\times \Bbb R^+\eqno (1.17)$$ so that $$dim\,\widetilde X
= dim\,X +2\eqno (1.18)$$ Since we are dealing here with a direct product of manifolds, functions, forms
and vector fields on $L^1$ and $\Bbb R^+$ have obvious extensions to $\widetilde X$ and we will freely
 use the same notation for the extensions. Let $r\in C^{\infty}(\Bbb R^+)$ be the natural
coordinate function. That is, $r(t) = t$ for any $t\in \Bbb R^+$. Let $\omega_{\widetilde X}$ be the
real closed (in fact exact) 2-form on $\widetilde X$ defined by putting $\omega_{\widetilde X}=
d(r\,\alpha)$ so that
$$\omega_{\widetilde X} = dr\wedge \alpha + r\,\widetilde {\omega}\eqno (1.19)$$ where we
have put $\widetilde {\omega} = \tau^*(\omega)$. 

By decomposing a tangent vector to $\widetilde X$ as a
sum of a tangent vector to $\Bbb R^+$ and a tangent vector to $L^1$ and then decomposing the latter into
a sum of a horizontal tangent vector (i.e., orthogonal to $\alpha$) and a vertical tangent vector (tangent
to a fiber of (1.5)) it follows easily that $\widetilde {\omega}$ is non-singular so one has \vs {\bf
Proposition 1.3.} {\it $(\widetilde X, \omega_{\widetilde X})$ is a symplectic manifold.}\vs It follows
immediately from (1.10) and (1.19) that the interior product of $\widetilde {\omega}$ by $\zeta$ equals
$dr$. Consequently one has \vs {\bf Proposition 1.4.} {\it The vector field $\zeta$ (see (1.10))
is Hamiltonian on $\widetilde X$. In fact $$\zeta = \xi_{r}\eqno (1.20)$$} \vs Let $\varphi\in
C^{\infty}(X)$ and $s\in S$. Note then, as a consequence of (1.11) and (1.12), Proposition 1.14 yields
the following Poisson bracket relations in $C^{\infty}(\widetilde X)$, $$\eqalign{[r,\widetilde
\varphi] &= 0\cr [r,\widetilde s]&= 2\,\pi\,i\,\, \widetilde s\cr}\eqno (1.21)$$
It is obvious from the sentence after (1.10) that the Hamiltonian flow of
$\zeta$ is just the free action of
$U(1)$ on
$\widetilde X$ defined by its principal bundle action on $L^1$. As a consequence of (1.19) one then
has \vs {\bf Proposition 1.5.} {\it Marsden-Weinstein reduction of $(\widetilde X,\omega_{\widetilde X})$
by the function $r$ at the value $r =1$ is isomorphic to the original symplectic
manifold $(X,\omega)$.} \vs Since $(X,\omega)$ arises from $(\widetilde X,\omega_{\widetilde X})$ by
symplectic reduction we can speak of the above construction of $(\widetilde X,\omega_{\widetilde X})$
from $(X,\omega)$ as symplectic induction. \vs 1.4. If $\eta$ is a vector field on a manifold then
$\iota(\eta)$ will denote the operator of interior product. We continue with the assumptions of
\S1.3.\vs {\bf Proposition 1.6.} {\it Let $\varphi\in C^{\infty}(X)$. Then $$\xi_{\widetilde {\varphi}} =
{1\over r}\,\widetilde {\xi_{\varphi}}\eqno (1.22)$$ so that if $\psi\in C^{\infty}(X)$ then
$$[\widetilde \varphi,\widetilde \psi] = {1\over r} \widetilde {[\varphi,\psi]}\eqno (1.23)$$ }\vs
{\bf Proof.} Since
$\widetilde {\xi_{\varphi}}$ is horizontal in $L^1$ clearly $\iota({1\over
r}\,\widetilde {\xi_{\varphi}})$ annihilates $d\,r\wedge \alpha$. But, from the
definition of $\xi_{\varphi}$, obviously
$$\iota(\widetilde {\xi_{\varphi}})\,\widetilde {\omega} = d\widetilde \varphi\eqno (1.24)$$ Thus
$$\eqalign{\iota({1\over r}\,\widetilde {\xi_{\varphi}})\omega_{\widetilde X}&= {1\over
r}\,\iota(\widetilde {\xi_{\varphi}})r\widetilde \omega\cr &= d\widetilde \varphi\cr}$$ But
this establishes (1.22). The equality (1.24) is immediate from the definition of Poisson bracket. QED\vs
The main part, (1.27), of the following result says that prequantization in
$X$ is ordinary Poisson bracket in the induced symplectic manifold $\widetilde X$. \vs {\bf Theorem 1.7.}
{\it The map $$C^{\infty}(X)\to C^{\infty}(\widetilde X),\qquad \varphi\mapsto r \widetilde
\varphi\eqno (1.25)$$ is a monomorphism of Poisson Lie algebras. Moreover $\eta_{\varphi}$ (see
Theorem 1.2) is a Hamiltonian vector field on $\widetilde X$ for any $\varphi\in C^{\infty}(X)$. In fact
$$ \eta_{\varphi} = \xi_{r \widetilde \varphi} \eqno (1.26) $$ Finally for any $s\in S$ and
$\varphi\in C^{\infty}(X)$ one has
$$\widetilde {\pi(\varphi)(s)} = [r\,\widetilde \varphi, \widetilde s]\eqno (1.27)$$}\vs {\bf
Proof.} Let $\varphi,\psi\in C^{\infty}(X)$. Then, by (1.21), $$[r\,\widetilde \varphi,
r\,\widetilde \psi] = r^2\,[\widetilde \varphi,\widetilde \psi]$$ But then
$[r\,\widetilde \varphi, r\,\widetilde \psi]= r \widetilde {[\varphi,\psi]}$ by (1.23).
This proves the first statement of the theorem. But now by (1.13) and (1.15)
one has $$\eqalign{\iota(\eta_{\varphi})\omega_{\widetilde X}&= \iota(\widetilde {\xi_{\varphi}} +
\widetilde \varphi\,\zeta)(dr\wedge\alpha + r\,\widetilde \omega)\cr &= \iota(\widetilde
{\xi_{\varphi}})r \widetilde \omega + \iota(\widetilde \varphi\,\zeta)(dr\wedge\alpha)\cr &=
r\,d\widetilde \varphi + \widetilde \varphi\,dr\cr}\eqno (1.28)$$ by (1.24) and (1.10).
But the right-hand side of the last line of (1.28) is just $d(r\,\widetilde \varphi)$. This proves
(1.26). But then (1.27) follows from (1.16). QED \vs {\bf Remark 1.8.} The construction of the
induced symplectic manifold $(\widetilde X,\omega_{\widetilde X})$ of course depends upon the
choice of the connection $\nabla$. The set (with an obvious equivalence relation) of all such
connections is a principal homogeneous space for the character group $\widehat {\pi_1(X)}$ of the
fundamental group $\pi_1(X)$ of $X$. See Theorem 2.5.1 in [K-1]. In particular the connection, and hence
$(\widetilde X,\omega_{\widetilde X})$, is unique if $X$ is simply-connected.
\vskip 1.5pc
\centerline{\bf 2. The coadjoint orbit case and a minimal Kostant-Sekiguchi correspondence}\vskip 1.5pc
2.1. We now consider the case of
\S1 where
$X$ is a coadjoint orbit of a connected Lie group $K$ and $\omega$ is the KKS symplectic form. We do not
assume now that
$K$ is compact but in the main application $K$ will be compact. Let
$\k = Lie\,K$ and let
$\k^*$ be the dual space to $\k$. Let
$\nu\in
\k^*$ and let $X$ be the $K$-coadjoint orbit of $\nu$ so that as $K$-homogenous spaces $$X =
K/K_{\nu}\eqno (2.1)$$ where $K_{\nu}$ is the isotropy subgroup at $\nu$. The action of $k\in K$ on
$\mu\in X$ is denoted by $k\cdot \mu$ and one has $\langle k\cdot \mu, y\rangle = \langle \mu,Ad\,k^{-1}
(y)\rangle$ for any $y\in \k$. For any
$x\in \k$ let
$\xi^x$ be the v.f. on $X$ defined so that for any $f\in C^{\infty}(X)$ and $\mu\in X$ one has
$(\xi^x\,f)(\mu) = {d\over dt}((exp(-tx)\cdot\mu)|_{t=0}$. Then $x\mapsto \xi^x$ defines the
infinitesimal action of $\k$ on $X$ corresponding to the group action of $K$ on $X$. The KKS symplectic
form $\omega$ on $X$ is such that for $x,y\in \k$ and $\mu\in X$ then $$\omega(\xi^x,\xi^y)(\mu) =
\langle\mu,[y,x]\rangle\eqno (2.2)$$ The map $$\k\to C^{\infty}(X),\qquad x\mapsto \varphi^x\eqno (2.3)$$
is a homomorphism of Lie algebras (using Poisson bracket in $C^{\infty}(X)$) where, for $x\in \k$, the
function
$\varphi^x$ on
$X$ is defined by
$$\varphi^x(\mu) = \langle \mu,x\rangle\eqno (2.4)$$ for any $\mu \in \k^*$. Furthermore $\xi^x$ is a
Hamiltonian vector field on $X$, for any $x\in \k$, and, in fact (see (5.3.5) in [K-1]), $$\xi^x =
\xi_{\varphi^x}\eqno (2.5)$$ Now $i\,\Bbb R = Lie\,U(1)$ and $$2\,\pi\,i\,\,\nu:\k_{\nu} \to i\,\Bbb
R\eqno (2.6)$$ is a homomorphism of Lie algebras where $\k_{\nu} = Lie\,K_{\nu}$. Note that $K_{\nu}$ may
not be connected. By Corollary 1 to Theorem 5.7.1 in [K-1] one has $[\omega]\in H^2(X,\Bbb R)$ is integral
if  there exists a character $$\chi:K_{\nu}\to U(1)\eqno (2.7)$$ whose differental is given by
$$d\,\chi = 2\,\pi\,i\,\,\nu|\k_{\nu}\eqno (2.8)$$ \vskip .5pc {\bf Remark 2.1.} In Corollary 1 to
Theorem 5.7.1 in [K-1] it is assumed that $K$ is simply-connected. In such a case the
statement of the corollary is that $[\omega]\in H^2(X,\Bbb R)$ is integral if and only if there
exists $\chi$ satisfying (2.7) and (2.8). The reference to the corollary is valid in the general
case under consideration here since (2.7) and (2.8) for
$K$ obviously imply these conditions for the simply-connected covering group of $K$. Furthermore, by the
corollary, the choice of $\chi$ determines $L^1$ and the connection 1-form $\alpha$ in $L^1$. See
(2.9) below. \vs Henceforth we assume (2.7) and (2.8) so
that the assumptions of
\S 1 are satisfied. Noting that
$K$ is a principal
$K_{\nu}$-bundle over
$X$ one constructs 
$L^1$ as the associated bundle given by $$L^1 = K\times_{K_{\nu}} U(1)\eqno (2.9)$$ where the action
of 
$K_{\nu}$ on $U(1)$ is given by the homomorphism $\chi$. See \S 5.7 in [K-1] where we have replaced $\Bbb
C^*$ by $U(1)$. 

By (2.9) the space $L^1$ is clearly homogeneous for the product group $K\times U(1)$ (using left
translation for $K$). In fact as homogeneous spaces one has an isomorphism $$(K\times U(1))/H\cong L^1$$
where $H=\{(k,\chi(k)^{-1})\mid k\in K_{\nu}\}$. In particular this gives rise to a $K\times
U(1)$-surjection
$\sigma:K\times U(1) \to L^1$. The connection 1-form $\alpha$ on $L^1$ may then be given by the
equation $$\sigma^*(\alpha) = (\alpha_{\nu},{d\theta\over 2\,\pi})\eqno (2.10)$$ where $\alpha_{\nu}$ is
the left invariant 1-form on $K$ whose value at the identity is $\nu$. See (5.7.3) and (5.7.4) in
[K-1]. \vs 2.2. We are assuming that $X$ is a coadjoint orbit of a connected
Lie group $K$ so that one has a Lie algebra homomorphism (2.3). We are assuming (2.7) and (2.8) so
that the induced symplectic manifold $(\widetilde X,\omega_{\widetilde X})$ exists. We now consider the
possibility that $(\widetilde X,\omega_{\widetilde X})$ may be the coadjoint orbit of some larger Lie
group
$G$. In such a case one would have a Lie algebra homomorphism $$\g\to C^{\infty}(\widetilde X)\eqno
(2.11)$$ where $\g = Lie\,G$. 

Assume now that $K$ is a compact connected Lie group. Let
$K_{\Bbb C}$ be a complex reductive Lie group having $K$ as a maximal compact subgroup so that
we can take $\k_{\Bbb C}= Lie\,K_{\Bbb C}$ where $\k_{\Bbb C}$ is the complexification of $\k$. Let
$B_{\k}$ be a non-singular symmetric $K_{\Bbb C}$-invariant bilinear form on $\k_{\Bbb C}$ which
is negative definite on $\k$. We may regard the complexification $\k^*_{\Bbb C}$ of $\k^*$ as the dual to
$\k_{\Bbb C}$. The bilinear form $B_{\k}$ defines a $K_{\Bbb C}$-linear isomorphism $$\gamma: \k^*_{\Bbb
C}\to
\k_{\Bbb C}\eqno (2.12)$$ Obviously $\gamma(\k^*) = \k$. It follows then that $K_{\nu}$ is just the
stabilizer of $\gamma(\nu)$ in $K$ under the adjoint representation so that $K_{\nu}$ is connected and,
as one knows,
$X$ is simply-connected. In addition one knows that $K_{\nu}$ contains a maximal torus $T$ of
$K$. One has $$\gamma(\nu)\in \tt\eqno (2.13)$$ where $\tt = Lie\,T$ since $\gamma(\nu)$ is obviously
central in $\k_{\nu}$. 

Let $\hh = i\tt$. The complexification $\hh_{\Bbb C} = \tt_{\Bbb C}$ is a Cartan subalgebra of $\k_{\Bbb
C}$ and, identifying the real dual $\hh^*$ of $\hh$ with $\gamma^{-1}(\hh)$ one has $\Lambda\s \hh^*$
where
$\Lambda$ is the $T$-weight lattice. But since
$2\,\pi\,i\,\,\nu|\tt$ exponentiates to a character of $T$, by (2.7), one has that $\lambda\in \Lambda$,
by (2.13), where
$$\lambda= 2\,\pi\,i\,\,\,\nu\eqno (2.14)$$ 

Let $\pi_{\lambda}:K_{\Bbb C} \to Aut\,V_{\lambda}$ be the irreducible representation of $K_{\Bbb C}$
with extremal $\hh_{\Bbb C}$-weight $\lambda$ (i.e., any $\hh_{\Bbb C}$-weight of $\pi_{\lambda}$ lies in
the convex hull of the Weyl group orbit of $\lambda$). Obviously the complexification $(\k_{\nu})_{\Bbb
C}$ is the centralizer of $\gamma(\lambda)$ in $\k_{\Bbb C}$. Furthermore $\gamma(\lambda)\in
\hh$ and hence
$\gamma(\lambda)$ is a hyperbolic element of $\k_{\Bbb C}$. \vs {\bf Remark 2.2.} To any 
hyperbolic element $y$ of a complex reductive Lie algebra $\ss$ one associates a
parabolic Lie subalgebra $q_y(\ss)$ of $\ss$ characterized as follows: the 
centralizer
$\ss^y$ of
$y$ in
$\ss$ is a Levi factor of $q_y(\ss)$ and the nilradical of $q_y(\ss)$ is the span of the eigenvectors of
$ad\,y$ belonging to positive eigenvalues. \vs In the notation of Remark 2.2
put $\q_{\nu} = \q_{\gamma(\lambda)}(\k_{\Bbb C})$ and let 
$Q_{\nu}\s K_{\Bbb C}$ be the parabolic subgroup corresponding to $\q_{\nu}$. One readily has $$K_{\nu} =
K\cap Q_{\nu}\eqno (2.15)$$ so that, since $X\cong K/K_{\nu}$ as $K$-homogeneous spaces the action of
$K$ on $K_{\Bbb C}/Q_{\nu}$ induces a $K$-isomorphism $$X\to K_{\Bbb C}/Q_{\nu}\eqno (2.16)$$

A non-zero vector in a finite dimensional complex irreducible $K$ (and hence $K_{\Bbb C}$)-module is
called an extremal weight vector if it is a weight vector for an extremal weight of some Cartan
subalgebra of $\k_{\Bbb C}$. The Cartan subalgebra can always be chosen so that it is the
complexification of a Cartan subalgebra of $K$. One knows that the set of all extremal weight vectors is
of the form
$K_{\Bbb C}\,\Bbb C^* \,v$ and in fact is of the form $K\,\Bbb C^*\,v$ where $v$ is an extremal weight
vector. 

Let
$E\s V_{\lambda}$ be the variety of extremal weight vectors in
$V_{\lambda}$. Thus if
$0\neq v_{\lambda}$ is a weight vector for the $\hh_{\Bbb C}$ extremal weight $\lambda$ then $E=
\pi_{\lambda}(K_{\Bbb C})\,\Bbb C^*\,v_{\lambda}$. Let
$Proj(V_{\lambda})$ be the projective space of $V_{\lambda}$ and let $Proj(E)$ be the subvariety of
$Proj(V_{\lambda})$ defined by $E$. One knows that $Proj(E)$ is the unique closed $K_{\Bbb C}$ orbit in 
$Proj(V_{\lambda})$. Furthermore if $p_{\lambda}\in Proj(E)$ is the point corresponding to $\Bbb
C^*\,v_{\lambda}$ then $Q_{\nu}$ is the isotropy group at $p_{\lambda}$ and hence (Borel-Weil theory)
the isomorphism (2.16) defines a $K$-isomorphism $$X\to Proj(E)\eqno (2.17)$$ 

Let $\{u,v\}$ be a Hilbert space structure ${\cal H}_{\lambda}$ in $V_{\lambda}$ which is invariant
under the action of $\pi_{\lambda}(K)$ and let $E^1 =\{v\in V_{\lambda}\mid \{v,v\} =1 \}$. We may
choose $v_{\lambda}$ so that $v_{\lambda}\in E^1$. Now if $\chi$ is defined (uniquely since $K_{\nu}$ is
connected) as in (2.7) and (2.8) so $L^1$ is given by (2.9) let $$\phi:L^1\to E^1\eqno (2.18)$$ be
the $K\times U(1)$-map given by $$\phi(k,c) = c\,\pi_{\lambda}(k)(v_{\lambda})\eqno (2.19)$$ where
$(k,c)\in K\times U(1)$. Furthermore if $t\in \Bbb R^+$ let $$\widetilde {\phi}:\widetilde {X}\to E\eqno
(2.20)$$ be the $K\times U(1)\times \Bbb R^+$ map given by $$\widetilde {\phi}(q,t) = t\,
\phi(q)\eqno (2.21)$$ where $(q,t)\in L^1\times \Bbb R^+$. The following is an immediate consequence of
(2.17).\vs {\bf Proposition 2.3.} {\it The maps $\phi$ and $\widetilde {\phi}$ (see (2.18) and (2.20))
are, respectively,
$K\times U(1)$ and $K\times U(1)\times \Bbb R^+$ isomorphisms.}\vs {\bf Remark 2.4.} Note that if
$\nu\neq 0$ then $K_{\Bbb C}$ operates transitively on $E$. This is clear since if $(K_{\nu})_{\Bbb
C}$ is the subgroup of $K_{\Bbb C}$ corresponding to $(\k_{\nu})_{\Bbb C}$ then obviously
$(K_{\nu})_{\Bbb C}$ operates transitively on $\Bbb C^*\,v_{\lambda}$. \vskip 1pc 2.3. We now assume that
$\g$ is a real simple Lie algebra of non-compact type and $G$ is a corresponding Lie group with finite
center. We now choose
$K$ and $\k$ of \S 2.2 so that $K$ is a maximal compact subgroup of $G$. Consequently there is a space
$\p\s
\g$ of hyperbolic elements such that $$\g = \k +\p\eqno (2.22)$$ is a Cartan decomposition of $\g$. Let
$\a\s\p$ be a maximum abelian subalgebra. Let $\Delta\s \a^*$ be the set of restricted roots and for
each $\beta \s \Delta$ let $\g_{\beta}\s \g$ be the corresponding restricted root space. Let
$\Delta_+\s \Delta $ be a choice of a positive root system and let $\n = \sum_{\beta\in
\Delta_+}\g_{\beta}$ so that $$\g = \k + \a + \n\eqno (2.23)$$ is an Iwasawa decomposition of $\g$.
Let $\n_- = \sum_{\beta\in \Delta_+}\g_{-\beta}$ and let $\m$ be the centralizer of $\a$ in $\k$ so that
one has the direct sum $$\g = \m + \a + \n +\n_-\eqno (2.24)$$ If
$\beta,\,\phi\in \Delta$ we will put $\beta \geq \phi$ if $\beta-\phi$ is a sum of elements in
$\Delta_+$. Obviously $\beta \geq \phi$ and $\phi \geq \beta$ if and only if $\beta = \phi$.
Let $\psi\in \Delta_+$ be a maximal element with respect to this ordering. Obviously $\g_{\psi}\s
Cent\,\n$. But since 
$Cent\,\n$ is clearly stable under the action of $ad\,\a$ it follows that $Cent\,\n$ is spanned by 
$\g_{\beta}' = \g_{\beta}\cap Cent\,\n$, over all $\beta\in \Delta_+$.
\vs {\bf Proposition 2.5.} {\it The restricted root $\psi$ is the unique maximal element in $\Delta_+ $
and 
$$Cent\,\n = \g_{\psi}$$ so that, in particular, $Cent\,\n$ is a restricted root space. } \vs
{\bf Proof.} Assume
$\g_{\beta}'\neq 0$ for some $\beta\in \Delta_+$. Since $\m$ normalizes $\n$ and commutes with $\a$ it
follows that $\g_{\beta}'$ is stable under the adjoint action of $\m + \a + \n$. Thus, corresponding to
the adjoint action, for the real enveloping algebras one has
$U_{\Bbb R}(\n_-)\,\g_{\beta}'=  U_{\Bbb R}(\g)\,\g_{\beta}'$ by (2.24), and the PBW theorem. But
by simplicity one has $U_{\Bbb R}(\g)\,\g_{\beta}' = \g$. Thus
$\beta\geq \psi$. But similarly $\psi\geq \beta$. Hence $\psi
=\beta$. QED\vs Let $W_{\a}$ be the restricted Weyl group operating in $\a$ and let $s_{\psi}\in W_{\a}$
be the reflection defined by $\psi$. Thus $s_{\psi}$ is the identity on the hyperplane $\a_{\psi}=
\{x\in \x\mid \psi(x) = 0\}$. Let $x_{\psi}$ be the unique element in $\a$ such that $s_{\psi}\,x_{\psi}
= -  x_{\psi}$ and $$\psi(x_{\psi})= 2\eqno (2.25)$$ 

The Killing form on the complexification $\g_{\Bbb
C}$ of $\g$ is real on $\g$ and positive definite on $\p$. Let $B$ be the positive multiple $(x,y)$ of
the Killing form normalized by the condition that $$(x_{\psi},x _{\psi}) = 2\eqno (2.26)$$ Let $B_{\k}$ of
\S 2.2 be defined so that $B_{\k} = B|\k$. Let $\theta$ be the complex Cartan involution on $\g_{\Bbb C}$
corresponding to the complexified Cartan decomposition $\g_{\Bbb C} = \k_{\Bbb C} + \p_{\Bbb C}$. Let
$\sigma$ and $\sigma_u$ be, respectively,  the conjugate linear involutions of $\g_{\Bbb C}$ defined
by the real forms $\g$ and $\g_u = k + i\,\p$. One notes that $\g_u$ is a compact form of $\g_{\Bbb C}$
so that $B$ is negative definite on $\g_u$. Consequently ${\cal H}$ is an $Ad\,\g_u$-invariant Hilbert
space structure $\{x,y\}$ on $\g_{\Bbb C}$ where $$\{x,y\} = - (x,\sigma_u\,y)\eqno (2.27)$$ It is
immediate that $\sigma$ commutes with $\theta$ and that $\sigma_u = \theta\,\sigma\,$. In particular the
restriction of ${\cal H}|\g$ defines a $K$-invariant real Hilbert space structure on $\g$ and that for
$x,y\in \g$, $$\{x,y\} = -(x,\theta y)\eqno (2.28)$$ Since $\theta$ is minus the identity on $\a$ one
immdiately has $$\theta (\g_{\psi}) = \g_{-\psi}\eqno (2.29)$$ Let $e\in \g_{\psi}$ be such that
$\{e,e\} = 1$. That is, $$(e,\theta e) =-1\eqno (2.30)$$ by (2.28).\vs {\bf Proposition 2.6.} {\it
$(x_{\psi},e,-\theta\, e)$ is an S-triple.} \vs {\bf Proof.} One has $[x_{\psi},e] = 2\,e$ by (2.25). But
$-\theta e\in \g_{-\psi}$ by (2.29) so that $[x_{\psi},-\theta\, e] = -2\,(-\theta\, e)$. Let $y =
[e,-\theta\, e]$. Then $[y,\a] = 0$ by (2.29). On the other hand clearly $\theta(y) = -y$ so that
$y\in\p$. Since $\a$ is maximally commutative in $\p$ this implies that $y\in \a$. But if $x\in \a$ then
$$\eqalign{(x,[e,-\theta\, e])&= ([x,e],-\theta\,e)\cr &=\psi(x)(e,-\theta\,e)\cr &=\psi(x)\cr}\eqno
(2.31)$$ by (2.30). This implies that $y$ is $B$-orthogonal to $\a_{\psi}$ so that $y=r\,x_{\psi}$ for
some $r\in \Bbb R$. But putting $x= x_{\psi}$ in (2.31) one has $(x_{\psi},y) = 2$ by (2.25). However
$(x_{\psi},x_{\psi}) = 2$ by (2.26). Thus $r=1$ so that $[e,-\theta\, e] = x_{\psi}$. QED\vs Let
$\u_{\psi}$ be the complex TDS spanned by $x_{\psi},e,$ and $-\theta\, e$ over $\Bbb C$. \vs 2.4. Let the
notation be as in
\S2.3. Let
$d = dim\,Cent\,\n$ so that $$d = dim\,\g_{\psi}\eqno (2.32)$$ by Proposition 2.5.\vs {\bf Proposition
2.7.} {\it The maximal eigenvalue of $ad\,x_{\psi}$ on $\g_{\Bbb C}$ is 2 and $(\g_{\psi})_{\Bbb C}$ is
the corresponding eigenspace. In particular the multiplicity of the eigenvalue 2 of $ad\,x_{\psi}$ is
$d$. Furthermore if $\beta \in \Delta-\{\psi\}$ then $ \beta(x_{\psi})\in \{0,1\}$ so that the spectrum
of $ad\,x_{\psi}$ in $\n_{\Bbb C}$ is non-negative and in fact the spectrum is contained in the set
$\{0,1,2\}$. The spectrum of $ad\,x_{\psi}$ on $\g_{\Bbb C}$ is contained in the set
$\{2,1,0,-1,-2\}$.}\vs {\bf Proof.} From the representation theory of a TDS (e.g. $\u_{\psi}$) one has
$\beta(x_{\psi})\in \Bbb Z$ for any $\beta\in \Delta_+$. On the other hand if $\beta\in \Delta_+$ then
since $\beta + \psi$ cannot be a restricted root one has $[e,\g_{\beta}] = 0$. Thus from the
representation theory of $\u_{\psi}$ one has $$\beta(x_{\psi})\in \Bbb Z_+\eqno (2.33)$$ But now
$\psi(x_{\psi}) = 2$ by (2.25) so that $(\g_{\psi})_{\Bbb C}$ is contained in the eigenspace of
$ad\,x_{\psi}$ for the eigenvalue 2.  But now using notation and the argument in the proof of Proposition
2.4 one has
$$\g_{\Bbb C} = (U_{\Bbb R}(\n_-)\,\g_{\psi})_{\Bbb C}\eqno (2.34)$$ so that 2 is the maximal
eigenvalue of $ad\,x_{\psi}$ by (2.33).  On the other hand if $\beta\in \Delta_+$ and
$[\g_{-\beta},\g_{\psi}]\neq 0$ we assert that $\beta(x_{\psi}) >0$. Note that the assertion
implies all the statements of the proposition, by (2.33) and (2.34). Assume the assertion is
false so that
$\beta(x_{\psi})= 0$. But then $\beta\neq \psi$ and hence $\psi-\beta\in \Delta$. But then
$\beta-\psi\in \Delta$ and $(\beta-\psi)(x_{-\psi}) = -2$. From the representation theory of a TDS one
has that $[e,[e,\g_{\beta-\psi}]]\neq 0$. But this implies that $\beta +\psi\in\Delta$ contradicting the
maximality of $\psi$. QED\vs Let $e_i,\,i=1,\ldots,d$, be an orthonormal basis of $\g_{\psi}$ with
respect to ${\cal H}|\g_{\psi}$. We assume that the basis is chosen so that $e_d = e$. Under the
adjoint action of $\u_{\psi}$ the element $e_i$ clearly generates a 3-dimensional irreducible
representation
$\u_i$ of $\u_{\psi}$ since $[e,e_i]= 0$. Of course $$\u_d = \u_{\psi}\eqno (2.35)$$\vskip .5pc {\bf
Remark 2.8.} Note that, as a consequence of Proposition 2.7, if $$\u = \sum_{i=1}^d \u_i\eqno (2.36)$$
then
$\u$ is the primary component in
$\g_{\Bbb C}$ for the 3-dimensional irreducible representation of $\u_{\psi}$ under the adjoint action
and that any irreducible component in $\g_{\Bbb C}/\u$ has dimension 1 or 2.\vs The following lemma is
well known and is readily established using the commutation relations of an S-triple. \vs {\bf Lemma
2.9.} {\it Assume that $\v$ is a complex TDS and $(x',e',f')$ is an S-triple whose elements span $\v$.
Then $(h,v,w)$ is an S-triple also spanning $\v$ where $$\eqalign{h&= i(e'-f')\cr v&= 1/2(ix' + e' +
f')\cr w&= 1/2(-ix' +e' +f')\cr}\eqno (2.37) $$ Furthermore one recovers $(x',e',f')$ from $(h,v,w)$ by
$$\eqalign{x'&=-i(v-w)\cr e'&= 1/2(-ih + v + w)\cr f'&= 1/2(ih + v + w)\cr}\eqno (2.38)$$}\vs We
apply Lemma 2.9 for the case where $\v = \u_{\psi}$ and $(x',e',f') = (x_{\psi},e,-\theta\, e)$. Let
$$\eqalign{h&= i(e + \theta\,e)\cr v&=1/2(ix_{\psi} + e - \theta\,e)\cr w&=1/2(-ix_{\psi} +e -
\theta\,e)\cr}\eqno (2.39)$$ so that $(h,v,w)$ is an S-triple whose elements span $\u_{\psi}$. 

Obviously there exists an automorphism of $\u_{\psi}$ which carries $(x_{\psi},e,-\theta\, e)$ to
$(h,v,w)$. Since any automorphism of a complex TDS is inner it follows that $x_{\psi}$ and $h$ are
conjugate in $\g_{\Bbb C}$ by an element in $Ad\,(\g_{\psi})_{\Bbb C}$. In particular then, by
Proposition 2.7, the multiplicity of the eigenvalue 2 of $ad\,h$ in $\g_{\Bbb C}$ is $d$. \vs
{\bf Lemma 2.10.} {\it One has $h\in \k_{\Bbb C}$ so that we may write $$d = d_{\k} + d_{\p}\eqno
(2.40)$$ where $d_k$ is the multiplicity of the eigenvalue 2 of $ad\,h|\k_{\Bbb C}$ and $d_{\p}$ is the
multiplicity of the eigenvalue 2 of $ad\,h|\p_{\Bbb C}$. One has $d_{\p}\geq 1$. In fact $v,w\in \p_{\Bbb
C}$ and $[h,v] = 2\,v$.}\vs {\bf Proof.} Obviously $\theta\,h = h$ by (2.39) so that $h\in \k_{\Bbb C}$.
Similarly, $v,w\in \p_{\Bbb C}$ since $\theta\,v = -v$ and $\theta\,w= -w$ by (2.39), noting
that $x_{\psi}\in\a\s \p_{\Bbb C}$. One has 
 $[h,v] = 2\,v$ since $(h,v,w)$ is an S-triple. QED\vs We can strengthen Lemma 2.10. \vs {\bf Theorem
2.11.} {\it Let the notation be as in Lemma 2.10. Then $d_{\k} =  d-1$ and $d_{\p} =1$. That is, the
one dimensional subspace $\Bbb C v$ is the eigenspace of $ad\,h|\p_{\Bbb C}$ corresponding to the
eigenvalue $2$. Also 2 is the highest eigenvalue of $ad\,h|\p_{\Bbb C}$.}\vs {\bf Proof.} The 3
dimensional $\u_{\psi}$-modules $\u_i$ (see (2.36)) are of course equivalent to the adjoint
representation of $\u_{\psi}$. For $j=1,\ldots,d,$ let $\delta_j:\u_{\psi}\to \u_j$ be the
$\u_{\psi}$-equivalence normalized so that $\delta_je = e_j$. Note that $\delta_d$ is the
identity map. Let $v_j = \delta_jv$. It is then immediate from Proposition 2.7
that $\{v_j\}\,j=1,\ldots,d,$ is a basis of the $ad\,h$ eigenspace in $\g_{\Bbb C}$ for the eigenvalue
2. To prove the first statement of the theorem it suffices, by Lemma 2.10, to show that
$$v_j\in\k_{\Bbb C}\,\,\hbox{for}\,\,i=j,\dots,d-1\eqno (2.41)$$ But now by the S-triple
commutation relations $ix_{\psi} = i[\theta\,e,e]$, $-\theta\,e = -1/2\, [\theta\,e,x_{\psi}]$ and $e=
-1/2\,[e,x_{\psi}]$. Thus if $x_j,f_j\in\u_j$ are defined by $x_j = [\theta\,e,e_j] $ and $f_j=
-1/2[\theta\,e,x_j]$ one has $$v_j = 1/2\,(ix_j + e_j + f_j)\eqno (2.42)$$ by (2.39). On the other hand
$$\eqalign{[e,x_j]&= [e,[\theta\,e,e_j]]\cr &= [[e,\theta\,e],e_j]\,\,\hbox{(since $[e,e_j]=0$)}\cr &=
-[x_{\psi},e_j]\cr &=-2\,e_j\cr}\eqno (2.43)$$ Let $j\in \{1,\ldots,d-1\}$. To establish (2.41) we will
first prove that $x_j\in \k_{\Bbb C}$. In fact we will prove that $$x_j\in \m\eqno (2.44)$$ Since
$\theta|\a$ is minus the identity one has $\theta\,e\in \g_{-\psi}$. But then $x_j$ clearly commutes with
$\a$. But the centralizer of $\a$ in $\g$ is $\m + \a$. To prove (2.44) it obviously suffices to prove
 that $x_j$ is
$B$-orthogonal to $\a$. But $\{e_j,e\} = \{e_j,e_d\} = 0$. Thus $-(e_j,\sigma_u\,e) = -(e_j,\theta\,e)=
0$. But then if $y\in \a$ one has $(y,x_j)= (y,[\theta\,e,e_j])= ([e_j,y],\theta\,e)$. But $[e_j,y] =
-\psi(y)\,e_j$. Thus $(y,x_j)= 0$ establishing (2.44). To prove (2.41) it now suffices, by (2.42), to
prove that $\theta\,f_j = e_j$. But $\theta\,f_j= -1/2\,[e,x_j]$ since $\theta\,x_j = x_j$ by (2.44).
But then $\theta\,f_j = e_j$ by (2.43). This proves (2.41). Since $h$ and $x_{\psi}$ are conjugate the
final statement of Theorem 2.11 follows from Proposition 2.7. QED\vs 2.5. We retain the notation of \S
2.4 and we will the apply the results of \S 2.4 to the symplectic considerations of \S 2.2.\vs {\bf 
Proposition 2.12.} {\it One has $$(v,w) = 1\eqno (2.45)$$ Furthermore $w = -\sigma_u v$ so that
$$\{v,v\} =1\eqno (2.46)$$}\vs {\bf Proof.} Since $h$ and $x_{\psi}$ are conjugate one has $$(h,h) =
2\eqno (2.47)$$ by (2.26). But since $(h,v,w)$ is an S-triple one has $(v,w)= 1/2([h,v],w)$. But
$([h,v],w)= (h,[v,w]) = (h,h)$. Thus (2.47) implies (2.45). But now $v=1/2(ix_{\psi} + e -
\theta\,e)$ and $ w=1/2(-ix_{\psi} +e - \theta\,e)$ by (2.39). Recall $\sigma_u = \theta\,\sigma$. But
clearly $\sigma\,v = 1/2(-ix_{\psi} + e - \theta\,e)$ and hence $\theta\,\sigma\,v = 
1/2(ix_{\psi} - e + \theta\,e)$. Hence $-\sigma_u\,v = w$. But then (2.46) follows from (2.45) and
(2.27). QED\vs Of course $h$ is a hyperbolic element in $\k_{\Bbb C}$ and $\q_h(\k_{\Bbb C})$ is the
parabolic subalgebra of $\k_{\Bbb C}$ defined by $h$. See Remark 2.2. \vs {\bf Theorem 2.13.} {\it
Under the adjoint action of $\k_{\Bbb C}$ on $\p_{\Bbb C}$ the one dimensional subspace $\Bbb C v$ is
stable under $\q_h(\k_{\Bbb C})$. In fact for any $x \in \q_h(\k_{\Bbb C})$ one has $$[x,v] =
(h,x)\,v\eqno (2.48)$$}\vs {\bf Proof.} If $x$ is contained in the nilradical of $\q_h(\k_{\Bbb C})$ then
$[x,v]= 0$ by the last line in Theorem 2.11. On the other hand if $x\in (\k_{\Bbb C})^h$ (a Levi factor
of 
$\q_h(\k_{\Bbb C}))$ then $\Bbb C\, v$ is stable under 
$ad\,x$ by the multiplicity one statement in Theorem 2.11 of the eigenvalue 2 of $ad\,h$ in $\p_{\Bbb
C}$. This proves the first statement of Theorem 2.13. For $x\in \q_h(\k_{\Bbb C})$ let $f$ be the linear
functional on
$\q_h(\k_{\Bbb C}) $ defined so that $[x,v] = f(x)\,v$. But then $f(x) = ([x,v],w)$ by (2.45). However
$([x,v],w) = (x,[v,w]) = (x,h)$. Thus $f(x) = (h,x)$. QED \vs Let $z = -ih$ so that $z\in \k$ and $z = e +
\theta e$ by (2.39). Let $T\s K$ be a maximal torus such that $z\in \tt$. Then $h\in \hh$ where, as in
\S2.2, $\hh = i\,\tt$. Let
$\gamma:\k_{\Bbb C}^*\to
\k_{\Bbb C}$ be as in (2.12) and let
$\lambda\in \hh^*$ (recalling the identification $\hh^* = \gamma^{-1}(\hh)$) be such that
$$\gamma(\lambda) = h\eqno (2.49)$$ Recall that $\Lambda\s \hh^*$ is the $T$-weight lattice. 
Note that $\hh_{\Bbb C} \s (\k_{\Bbb C})^h$ so that $\hh_{\Bbb C}\s \q_h(\k_{\Bbb C})$.
As an immediate consequence of
Theorem 2.13 one has
\vs {\bf Proposition 2.14.} {\it One has $\lambda\in \Lambda$. Furthermore the $K_{\Bbb C}$-module $V$
generated by $v$ with respect to the adjoint action of $K_{\Bbb C}$ on $\p_{\Bbb C}$ is irreducible and
is equivalent to $V_{\lambda}$ with $v$ corresponding to $v_{\lambda}$.}\vs Henceforth we will
identify
$V$ with $V_{\lambda}$ and $v$ with $v_{\lambda}$. One has $\pi_{\lambda}(k) y = Ad\,k(y)$ where $y\in
V$ and $k\in K_{\Bbb C}$. 

 Let $\nu = \lambda/2\,\pi\,i$ so that $\nu\in \k^* $. Let, as in \S2.2, $X$ be
the $K$-coadjoint orbit of $\nu$ so that $(X,\omega)$ is a symplectic $K$-homogeneous space where $\omega$
is the KKS symplectic form. One readily notes that $$Lie\,K_{\nu}=\k\cap (\k_{\Bbb C})^h\eqno (2.50)$$
so that $\k\cap (\k_{\Bbb C})^h$ is a compact form of the Levi factor $(\k_{\Bbb C})^h$ of
$\q_h(\k_{\Bbb C})$ where
$K_{\nu}$ is the isotropy group at $\nu$. Furthermore one knows (as a general fact about coadjoint orbits
of compact connected Lie groups) that
$K_{\nu}$ is connected so that (2.7) and (2.8) are satisfied where $\chi$ is the character on $K_{\nu}$
defined by the action of $K_{\nu}$ on $\Bbb C\,v_{\lambda}$. Thus, as in \S2.2, we can construct the
induced symplectic manifold $(\widetilde {X},\omega_{\widetilde {X}})$. 

As in \S 2.2 let $E\s V_{\lambda}$ be the variety of extremal weight vectors so that $E =
\pi_{\lambda}(K_{\Bbb C})\,\Bbb C^*\,v_{\lambda}$. Let ${\cal H}_{\lambda}$ be the $K$-invariant
Hilbert space structure in $V_{\lambda}$ given the restriction ${\cal H}|V_{\lambda}$. As in \S2.2,
$E^1$ is the space of vectors in $E$ having length 1 with respect to ${\cal H}_{\lambda}$. Note that
$$v_{\lambda}\in E^1\eqno (2.51)$$ by (2.46). We recall that Proposition 2.3 sets up a $K\times U(1)$
isomorphism $$L^1\to E^1\eqno(2.52)$$ and a $K\times U(1)\times \Bbb R^+$ isomorphism $$\widetilde
{X}\to E\eqno (2.53)$$ Let $K_{v}$ be the isotropy group at $v = v_{\lambda}$ for the action of $K$ on
$E^1$ and let $K_e$ be the isotropy group at $e\in \g$ for the adjoint action of $K$ on $\g$. \vs {\bf
Theorem 2.15.} {\it The following 3 subgroups of $K$ are equal. $$\eqalign{&(1)\,\,K_{v}\cr
&(2)\,\,K_{e}\cr &(3)\,\,\hbox{The centralizer of the TDS $\u_{\lambda}$ in $K$}\cr}$$}\vs {\bf Proof.
} Let $K'$ be the subgroup of $K$ defined by (3). Since $v,e\in \u_{\psi}$ one obviously has $K'\s K_v$
and
$K'\s K_e$. But clearly $Ad\,k$ commutes with $\theta$ for any $k\in K$. Thus $\theta (e)$ is fixed by
$Ad\,k$ for any
$k\in K_e$. But then all 3 elements of the S-triple in Proposition 2.6 are fixed by $Ad\,k$. Hence $K_e
= K'$ by the definition of $\u_{\psi}$ (see \S 2.3). But $\sigma_u$ also commutes with the adjoint action
of $K$ since $\k\s \g_u$. But $w = -\sigma_u\,v$ by Proposition 2.12. Thus $w$ is fixed by $Ad\,k$ for any
$k\in K_v$. Hence any such $k$ fixes the three elements of the S-triple $(h,v,w)$. But these elements
also span $\u_{\psi}$ . Thus $K' = K_v$. QED\vs Let $O$ be the $Ad\,G$ orbit of $e$ in $\g$ and let $O^1
= \{f\in O\mid \{f,f\} =1\}$. Note that $e\in O^1$ by the choice of $e$ in \S 2.3. Since $K$ operates
unitarily with respect to ${\cal H}$ the adjoint action of $k$ stabilizes $E^1$ and $O^1$. As a corollary
of Theorem 2.15 one has 
\vs  {\bf Theorem 2.16.} {\it The compact group $K$ operates transitively on $E^1$ and on $O^1$.
Furthermore these spaces are isomorphic as
$K$-homogeneous spaces. In fact $b$ is such an isomorphism where for any $k\in K$,
$$b(Ad\,k\,(v_{\lambda})) = Ad\,k\,(e)\eqno (2.54)$$}\vs {\bf Proof.} Recall $z = e + \theta\,e =
-ih\in \k$. But if $s\in \Bbb R$ one then has
$$\pi_{\lambda}(exp\,s\,z)v_{\lambda} = e^{-2\,s\,i}\,v_{\lambda} \eqno (2.55)$$ by (2.47) and (2.48).
In particular $k\mapsto \phi\,(k,1)$ in (2.19) surjects $k$ onto $E^1$ by Proposition 2.3. Consequently
$K$ operates transitively on $E^1$. On the other hand if $A$ and $N$ are the subgroups of $G$ which
correspond, respectively, to $\a$ and $\n$ then one has the group Iwasawa decomposition $G = KAN$. But
$e$ is fixed under the adjoint action of $N$ since $e\in Cent\,\n$ (see Proposition 2.5). On the other
hand
$Ad\,A\,(e) = \Bbb R^+\,e$. Thus $$O = Ad\,K\,\,\,\Bbb R^+\,e\eqno (2.56)$$ However clearly 
$\Bbb R^+\,e\cap O^1=\{e\}$. Hence $K$ operates transitively on $O^1$. The theorem then follows from the
equality of (1) and (2) in Theorem 2.15. QED \vs The variety $E$ is $K_{\Bbb C}$ homogeneous by
Remark 2.4. The $K_{\Bbb C}$ orbit $E$ in $\p_{\Bbb C}$ corresponds to the $G$-orbit $O$
in $\g$ by the Kostant-Sekiguchi correspondence (a correspondence of $K_{\Bbb C}$ nilpotent orbits in
$\p_{\Bbb C}$ and $G$ nilpotent orbits in $\g$. See [S]). Mich\`ele Vergne has proved the corresponding
orbits are $K$-diffeomorphisms (see [V]). The proof is highly non-trivial. However, by Theorem 2.16, in
the special case of 
$E$ and $O$ the diffeomorphism is transparent. Obviously by (2.21) and (2.56) one has
$K\times \Bbb R^+$-diffeomorphisms
$$\eqalign{&E^1\times \Bbb R^+\to E\,\,\,\qquad(f,t) \mapsto tf\cr &O^1\times \Bbb R^+\to
O\,\,\,\qquad(u,t)
\mapsto tu\cr}\eqno (2.57)$$ We may therefore extend the domain of definition of $b$
so that, using the notation of (2.57), one has a $K\times \Bbb R^+$-diffeomorphism $$b:E\to
O,\,\,\hbox{where}\,\,b\,(t\,f) = t\,b\,(f)\eqno (2.58)$$ \vskip .5pc

\rm
\centerline{\bf 3. The symplectic isomorphism $(\widetilde
X,\omega_{\widetilde X})\cong (Z,\omega_Z)$}\vskip 1.5pc
3.1. We continue with the notation of \S 2. Recalling (2.12) one
notes that since $B_{\k} = B|\k_{\Bbb C}$ and, since $\k_{\Bbb C}$
and $\p_{\Bbb  C}$ are $B$-orthogonal, the isomorphism $\g_{\Bbb
C}^*$ to $\g_{\Bbb C}$ defined by $B$ is an extension of (2.12).
The extension will also be denoted $\gamma$. Let $\varepsilon \in
\g^*$ be defined so that $$\gamma(\varepsilon) = e/\pi\eqno (3.1)$$
Let $Z\s \g^*$ be the $G$-coadjoint orbit of $\varepsilon$ and let
$\omega_Z$ be the KKS-symplectic form on $Z$. To avoid confusion
with the vector field $\xi^x$ on $X$ (see \S 2.1) defined by any
$x\in \k$, with respect to the coadjoint action of $K$ on $X$, we
will denote by
$\Xi^y$ the vector field on
$Z$ defined by any $y\in\g$ with respect to the coadjoint action
of $G$ on $Z$. The analogue of (2.2) is the formula
$$\omega_Z(\Xi^x,\Xi^y)(\rho) = \langle \rho,[y,x]\rangle 
\eqno (3.2)$$ for any $\rho\in Z$. In particular if $\rho = \varepsilon$ one has 
$$\omega_Z(\Xi^x,\Xi^y)(\varepsilon) = (1/\pi)\,(e,[y,x])\eqno (3.3)$$ 
Recall that $\omega$ is the KKS form on the $K$-coadjoint orbit $X$ of $\nu =
\lambda/2\,\pi\,i$. \vs
{\bf Lemma 3.1.} {\it Let
$x,y\in
\k$. Then
$$\omega_Z(\Xi^x,\Xi^y)(\varepsilon) =  \omega(\xi^x,\xi^y)(\nu)\eqno (3.4)$$}\vs {\bf
Proof.} $[y,x]$ is fixed by $\theta$ since $[y,x]\in\k$. Thus
$(1/\pi)\,(e,[y,x]) = ((e + \theta\,e)/2\,\pi,[y,x])$. But $e+\theta\,e =
(1/i)h$ by (2.39) and $\gamma(\lambda) = h$. But then $$\gamma(\nu) = (e +
\theta\,e)/2\,\pi\eqno (3.5)$$ Thus $(1/\pi)\,(e,[y,x]) = \langle \nu,[y,x]\rangle $.
But then (3.4) follows from (2.2) and (3.3). QED\vs Now recalling (2.9) the circle
bundle $L^1$ is given by $L^1 = K\times_{K_{\nu}}U(1)$ where the action of $K_{\nu}$
on $U(1)$ is given by the character (recall $K_{\nu}$ is connected)
$$\eqalign{\chi(exp\,x) &= e^{\lambda(x)}\cr &= e^{(h,x)}\cr}\eqno (3.6)$$ for any
$x\in \k_{\nu}$ (see (2.48)). Now by (1.17) one has $\widetilde {X} = L^1\times \Bbb
R^+$ so that $$\widetilde {X} = K\times_{K_{\nu}}U(1)\times \Bbb R^+\eqno (3.7)$$
Let $o\in \widetilde {X}$ be the point whose components are the identity in $K$, 1 in
$U(1)$ and 1 in $\Bbb R^+$ with respect to (3.7). Extend the domain of the bundle
projection
$\tau$ (see (1.6)) to $\widetilde {X}$ so that $\tau(q,t) = \tau(q)$ for $(q,t)\in
L^1\times \Bbb R^+$. One notes that $$\tau(o) = \nu\eqno (3.8)$$\vskip .5pc {\bf
Proposition 3.2.} {\it There exists a
$K\times
\Bbb R^+$ diffeomorphism
$$\beta:\widetilde {X}\to Z\eqno (3.9)$$ where $$\beta(t\,k\cdot o) = t \,(Coad\,\,k
(\varepsilon))\eqno (3.10)$$}\vs {\bf Proof.} Put $\beta = (1/\pi) \gamma^{-1}\circ
b\circ \widetilde {\phi}$. The result then follows from Proposition 2.3, Theorem
2.16, (2.58) and the invariance of
$B$. QED \vs 3.2. Our main objective will be to prove that $\beta:(\widetilde
{X},\omega_{\widetilde {X}})\to (Z,\omega_Z)$ is an isomorphism of symplectic
manifolds. 

Let $\beta_*:T(\widetilde {X})\to T(Z)$ be the diffeomorphism of tangent bundles
defined by the differential of $\beta$. Let $\beta_o$ be the restriction of
$\beta_*$ to the tangent space $T_o(\widetilde {X})$ so that
$$\beta_o:T_o(\widetilde {X})\to T_{\varepsilon}(Z)\eqno (3.11)$$ is a linear
isomorphism. For any $x\in \k$ let $\eta^x$ be the vector field on $\widetilde {X}$
defined by the action of $K$ on $\widetilde {X}$. Since $\beta$ is a $K$-map one has
$$\beta_*(\eta^x) = \Xi^x\eqno (3.12)$$ One can be very explicit about $\eta^x$.\vs
{\bf Proposition 3.3.} {\it Let $x\in \k$. Then using the notation of (1.13) and
(2.4) one has $$\eta^x = \widetilde {\xi^x} +\widetilde
{\varphi^x}\,\zeta\eqno (3.13)$$} \vs {\bf Proof.} As a vector field on $L^1$ (and
hence on $\widetilde {X}$) $\eta^x$ is characterized by the property that (1) it
commutes with the
$U(1)$-action, (2) it annihilates $\alpha$ by Lie differentiation and
(3) $$\tau_*(\eta^x) =\xi^x\eqno (3.14)$$ The result then follows from (1.13), (2.5)
and Theorem 1.2. QED\vs Let $\k^{\perp}_{\nu}$ be the $B$-orthocomplement of
$\k_{\nu}$ in $\k$ so that the map $$\k^{\perp}_{\nu}\to T_{\nu}(X)\eqno (3.15)$$
given by
$x\mapsto (\xi^x)_{\nu}$ is a linear isomorphism. Let $R_o$ be the space of
horizontal tangent vectors (i.e., orthogonal to $\alpha$) to $L^1$ at $o$ so that
$R_o$ has codimension 1 in
$T_o(L^1)$ and codimension 2 in $T_o(\widetilde {X})$. \vs {\bf Lemma 3.4.} {\it
One has $(\eta^x)_o\in R_o$ for any $x\in \k^{\perp}_{\nu}$ and the map
$$\k^{\perp}_{\nu} \to R_o,\,\qquad x\mapsto (\eta^x)_o\eqno (3.16)$$ is a linear
isomorphism.} \vs {\bf Proof.} Let $x\in k^{\perp}_{\nu}$. Then $\varphi^x(\nu) = 0$
by (2.4) since $\gamma(\nu)\in \k_{\nu}$. Thus $$(\eta^x)_o = \widetilde
{\xi^x}\eqno (3.17)$$ by (3.13). This proves that $(\eta^x)_o\in R_o$. But
since $\tau_*:R_o\to T_{\nu}(X)$ is clearly an isomorphism the remaining statements
of the proposition follow from (3.14) and the isomorphism (3.15).  QED\vs Let
$R_{\varepsilon} = \beta_o(R_o)$ so that (3.11) restricts to the linear isomorphism
$\beta_o:R_o\to R_{\varepsilon}$. Let $\omega_o$ be the symplectic
bilinear form on $T_o(\widetilde {X})$ induced by the retriction of
$\omega_{\widetilde {X}}$ to $T_o(\widetilde {X})$ and let $\omega_{\varepsilon}$ be
the symplectic bilinear form on $T_{\varepsilon}(Z)$ induced by the retriction of
$\omega_{Z}$ to $T_{\varepsilon}(Z)$. We wish to prove that
$$\beta_o:(T_o(X),\omega_o)\to (T_{\varepsilon}(Z),\omega_{\varepsilon})\eqno
(3.18)$$ is an isomorphism of symplectic vector spaces. We first establish \vs {\bf
Lemma 3.5.} {\it The restrictions $\omega_o|R_o$ and
$\omega_{\varepsilon}|R_{\varepsilon}$ are non-singular and
$$\beta_o:(R_o,\omega_o|R_o)\to
(R_{\varepsilon},\omega_{\varepsilon}|R_{\varepsilon})\eqno (3.19)$$ is an
isomorphism of symplectic vector subspaces.}\vs {\bf Proof.} Recalling the
definition of 
$\omega_{\widetilde {X}}$ (see (1.19)) it is clear that if $x,y\in \k^{\perp}_{\nu} $
then $\omega_{\widetilde {X}}(\eta^x,\eta^y)(o) = \widetilde {\omega}(\widetilde
{\xi^x},\widetilde {\xi^y})(o)$ by (3.17). But then $\omega_{\widetilde
{X}}(\eta^x,\eta^y)(o) = \omega(\xi^x,\xi^y)(\nu)$. Thus $\omega_o|R_o$ is
non-singular by the linear isomorphism (3.15). But then
$\omega_{\varepsilon}|R_{\varepsilon}$ is non-singular and (3.19) is an
isomorphism of symplectic vector subspaces by (3.4) and (3.12). QED\vs Now let
$R_o^{\perp}$ be the 2-dimensional orthocomplement of $R_o$ in $T_o(\widetilde
{X})$ with respect to $\omega_o$. It is clear from (1.19) and \S 1.2 that
$R_o^{\perp}$ is spanned by $\zeta_o$ and $(d/dr)_o$. It will be convenient for us
to modify this basis of $R_o^{\perp}$. Recall that $z\in \k$ is given by $z = e +
\theta\,e$.
\vs {\bf Lemma 3.6.} {\it One has $$(\eta^z)_o = -\zeta_o/\pi\eqno (3.20)$$ so
that
$(\eta^z)_o$ and $(-2\,r\,d/dr)_o$ are a basis of $R_o^{\perp}$. Furthermore
$$\omega_o((-2\,r\,d/dr)_o, (\eta^z)_o) = -2/\pi\eqno (3.21)$$} \vs {\bf Proof.} By
(2.49) one has $$\gamma(\nu) = z/2\,\pi\eqno (3.22)$$ so that $z\in \k_{\nu}$. Thus
$(\xi^z)_{\nu} = 0$. Hence $(\eta^z)_o = (\widetilde {\varphi^z}\zeta)_o$ by (3.13).
But $\widetilde {\varphi^z}(o) = \varphi^z(\nu)$ and $\varphi^z(\nu) =
(\gamma(\nu),z)$ by (2.4). Hence $\widetilde {\varphi^z}(o) = 1/2\pi\,(z,z)$ by
(3.22). But since $z = -ih$ one has $(z,z) = -2$ by (2.47). This proves (3.20). But
now by (1.19), since $r(o) = 1$, $$\eqalign{\omega_o((-2\,r\,d/dr)_o, (\eta^z)_o) &=
(dr\wedge\alpha) (-2\,r\,d/dr,\eta^z)(o)\cr &= -2/\pi\cr}$$ by (3.20) and (1.10).
This proves (3.21). QED \vs Let $\kappa$ be the Euler vector field on $Z$. Thus if
$f\in C^{\infty}(Z)$ and $\mu\in Z$ then $(\kappa\,f)(\mu) = d/dt\,f(\mu +
t\mu)|_{t=0}$. Clearly $$\beta_*(r\,\,d/dr) = \kappa\eqno (3.23)$$\vskip.5pc {\bf
Lemma 3.7.} {\it One has $$\beta_o((-2\,r\,\,d/dr)_o) =
(\Xi^{x_{\psi}})_{\varepsilon}\eqno (3.24)$$ (see Proposition 2.6).}\vs {\bf Proof.}
One has $[x_{\psi},e] = 2\,e$ by Proposition 2.6. But then
$$coad\,x_{\psi}(\varepsilon) = 2\,\varepsilon\eqno (3.25)$$ since
$\pi\,\gamma$ is an equivalence of $\g$-modules (see (3.1)). But then
$(\Xi^{x_{\psi}})_{\varepsilon} = -2\kappa_{\varepsilon}$. (See the beginning of \S
2.1 to explain the minus sign.) But then (3.24) follows from (3.23). QED \vs
{\bf Lemma 3.8.} {\it One has
$$\omega_{\varepsilon}((\Xi^{x_{\psi}})_{\varepsilon},(\Xi^z)_{\varepsilon}) =
-2/\pi\eqno (3.26)$$}\vs {\bf Proof.}
$$\eqalign{\omega_{\varepsilon}((\Xi^{x_{\psi}})_{\varepsilon},(\Xi^z)_{\varepsilon})
&= \omega_Z(\Xi^{x_{\psi}},\Xi^z)(\varepsilon)\cr &= \langle
\varepsilon,[z,x_{\psi}]\rangle\cr &= (1/\pi)\,(e,[z,x_{\psi}])\cr & = (1/\pi)\,([e,e
+\theta\,e],x_{\psi})\cr &= -(1/\pi)\,(x_{\psi},x_{\psi})\,\,\quad\hbox{by 
Proposition 2.6}\cr &= -2/\pi\,\,\quad\hbox{by (2.26)}\cr}$$ QED\vs We can now prove
\vs {\bf Theorem 3.9}. {\it The map $$\beta_o:(T_o(X),\omega_o)\to 
(T_{\varepsilon}(Z),\omega_{\varepsilon})\eqno
(3.27)$$ is an isomorphism of symplectic vector spaces.} \vs {\bf Proof.} 
Let $R_{\varepsilon}^{\perp}$ be the 2  dimensional orthocomplement of
$R_{\varepsilon}$ in $T_{\varepsilon}(Z)$ with respect to $\omega_{\varepsilon}$. We
assert (Assertion A) that $(\Xi^z)_{\varepsilon}$ and $(\Xi^{x_{\psi}})_{\varepsilon}$
is a basis of
$R_{\varepsilon}^{\perp}$. We first show that $(\Xi^z)_{\varepsilon}\in
R_{\varepsilon}^{\perp}$. To do this it suffices to show that
$$\omega_Z(\Xi^z,\Xi^x)(\varepsilon) = 0\eqno (3.28)$$ for all $x\in \k$ since
$R_{\varepsilon}$ is spanned by $(\Xi^x)_{\varepsilon}$ for $x\in \k_{\nu}^{\perp}$,
by Lemma 3.4 and Lemma 3.5. But if $x\in \k$ then
$$\eqalign{\omega_Z(\Xi^z,\Xi^x)(\varepsilon)&= 1/\pi(e,[x,z])\cr &= (1/2\,\pi) ((e +
\theta e),[x,z])\,\,\quad\hbox{since $[x,z]\in
\k$}\cr &= (1/2\,\pi)([z,z],x)\cr &=0\cr}$$ Thus $(\Xi^z)_{\varepsilon}\in
R_{\varepsilon}^{\perp}$.

Let $y\in \k_{\nu}^{\perp}$. Then
$$\eqalign{\omega_{Z}(\Xi^{x_{\psi}},\Xi^y)(\varepsilon)&= (1/\pi)(e,[y,x_{\psi}])\cr
&=(1/\pi)([x_{\psi},e],y)\cr &=(2/\pi)(e,y)\,\,\qquad\hbox{by Proposition 2.6}\cr
&=(1/\pi)(e + \theta\,e,y)\,\,\hbox{since $y\in \k$}\cr &=0\,\,\qquad\hbox{since
$z\in\k_{\nu}$ and $y\in \k_{\nu}^{\perp}$}\cr}$$ But then
$(\Xi^{x_{\psi}})_{\varepsilon}\in R_{\varepsilon}^{\perp}$ by Lemma 3.4 and
Lemma 3.5. This proves Assertion A since the left side of (3.26) is non-zero. But
now we assert (Assertion B) $\beta_o(R_o^{\perp})\s R_{\varepsilon}^{\perp}$ and
$$\beta_o:(R_o^{\perp},\omega_o|R_o^{\perp})\to
(R_{\varepsilon}^{\perp},\omega_{\varepsilon}|R_{\varepsilon}^{\perp})\eqno
(3.29)$$ is an isomorphism of symplectic subspaces. Indeed $\beta_o((\eta^z)_o) =
(\Xi^z)_{\varepsilon}$ by (3.12) and $\beta_o((-2\,r\,d/dr)_o) =
(\Xi^{z_{\psi}})_{\varepsilon}$ by (3.24). But then Assertion B follows from
Assertion A together with Lemma 3.6 and the fact that $-2/\pi$ occurs on the right
side of (3.21) and on the right side of (3.26).  But then the theorem follows from
the symplectic isomorphism (3.29) together with the symplectic isomorphism (3.19).
QED\vs The following is our main result.\vs {\bf Theorem 3.10}. {\it Let $G$ be
a connected non-compact Lie group with finite center such that $Lie\,\g$ is simple.
Let
$K$ be a maximal compact subgroup and let $X$ be the coadjoint orbit of $K$ 
defined as in \S 2.5. Let $(\widetilde {X},\omega_{\widetilde {X}})$ be the symplectic
manifold obtained from $X$ by symplectic induction so that $dim\,\widetilde {X} =
dim\,X +2$. See \S 1.3. Let $(Z,\omega_Z)$ be the coadjoint orbit of $G$, together
with its $KKS$ symplectic form, defined as in \S 3.1, and let $$\beta:(\widetilde
{X},\omega_{\widetilde {X}})\to (Z,\omega_Z)\eqno (3.30)$$ be the $K$-diffeomorphism
(special case of M. Vergne's theorem)  defined as in \S3.1 using the
Kostant-Sekiguchi correspondence. Then (3.30) is an isomorphism of symplectic
manifolds. } \vs {\bf Proof.} Note that $L^1$ (see (2.9)) is $K$-homogeneous
since (2.7) is surjective ($\nu\neq 0$ since $G$ is non-compact). Thus $\widetilde
{X}$ is $K\times \Bbb R^+$ homogeneous by (1.17). For any $(k,s)\in K\times \Bbb
R^+$ let $\ell_{k,s}$ be the diffeomorphism of $\widetilde {X}$ defined by the action
of $(k,s)$ on $\widetilde {X}$. Thus if $q\in \widetilde {X}$ then $\ell_{k,s}(q) =
s\,(k\cdot q)$. Let $m_{k,s}$ be the diffeomorphism of $Z$ defined so that if
$\mu\in Z$ then $m_{k,s}(\mu) = s\,(Coad\,k (\mu))$. It follows from (2.56) that
$O$ is $K\times \Bbb R^+$ homogeneous. Consequently $Z$ is $K\times \Bbb R^+$
homogeneous since (see (3.1)) $\pi\,\gamma:Z\to O$ obviously commutes with the action
of $K\times \Bbb R^+$. From the definition of $\beta$ (see (3.10)) it is immediate
that one has the commutative diagram equality of maps $$\beta \circ \ell_{k,s} = 
m_{k,s}\circ \beta\eqno (3.31)$$ Let $q\in \widetilde {X}$ and put
$p=\ell_{k,s}\,q$. We assert that the pullback
$$(\ell_{k,s})^*((\omega_{\widetilde {X}})_p) = s\,(\omega_{\widetilde {X}})_q\eqno
(3.32)$$ But recalling (1.19) it is obvious that $\alpha$ and $\widetilde
{\omega}$ are invariant under $\ell_{k,s}$. On the other hand clearly
$(\ell_{k,s})^*r = s\,r$ and hence $(\ell_{k,s})^*dr = s\,dr$. This proves (3.32). 

Now let $\mu\in Z$ and let $\rho = m_{k,s}(\mu)$. We assert that $$
(m_{k,s})^*((\omega_Z)_{\rho}) = s\,(\omega_{Z})_{\mu}\eqno (3.33)$$ Since
$\omega_{Z}$ is obviously $K$-invariant it suffices to prove (3.33) under the
assumption that $k$ is the identity of $K$. Making that assumption one has 
$\rho = s\,\mu$. But also 
$\Xi^y$, for any $y\in \g$, is invariant under $m_{k,s}$. Thus for any $y\in \g$ one
has
$$(m_{k,s})_*((\Xi^y)_{\mu}) = (\Xi^y)_{s\,\mu}\eqno (3.34)$$ But then for any
$x,y\in
\g$ one has $$\eqalign{(m_{k,s})^*((\omega_Z)_{s\,\mu})((\Xi^x)_{\mu},(\Xi^y)_{\mu})&
=(\omega_Z)_{s\,\mu}((\Xi^x)_{s\mu},(\Xi^y)_{s\mu})\cr &= \langle s\,\mu,
[y,x]\rangle\cr &= s\,(\omega_Z)_{\mu}((\Xi^x)_{\mu},(\Xi^y)_{\mu})\cr}$$ But this
proves (3.33). Let $p\in \widetilde {X}$ be arbitrary and let $\mu
=\beta\,p$. To prove the Theorem 3.10 it suffices to show that
$$\beta^*((\omega_Z)_{\mu}) = (\omega_{\widetilde {X}})_{p}\eqno (3.35)$$ By
transitivity there exists 
$(k,s)\in K\times \Bbb R^+$ such that $\ell_{k,s}p = o$. But then $m_{k,s}\mu =
\varepsilon$ by the commutativity equation (3.31) and in fact $$
\beta^*((\omega_Z)_{\mu}) =
(\ell_{k,s})^*((\beta_o)^*(((m_{k,s})^{-1})^*(\omega_Z)_{\mu}))$$ But
$(m_{k,s})^{-1})^*(\omega_Z)_{\mu} = s^{-1}\,(\omega_Z)_{\varepsilon}$ by (3.33).
But
$(\beta_o)^*(s^{-1}\,(\omega_Z)_{\varepsilon}) = s^{-1} (\omega_{\widetilde {X}})_o$
by Theorem 3.9. Finally $(\ell_{k,s})^*(s^{-1} (\omega_{\widetilde {X}})_o) =
(\omega_{\widetilde {X}})_{p}$ by (3.32). This proves (3.35). QED\vs 3.3. We wish to
characterize the varieties $Z$ and $X$ in more general terms.
Recall (E. Cartan's theory) that the non-compact symmetric space
$G/K$ is one of two types (1) non-Hermitian symmetric or (2) Hermitian symmetric. In
the non-Hermitian case
$\p_{\Bbb C}$ is $K_{\Bbb C}$-irreducible and $\k_{\Bbb C}$ is semisimple. In the
Hermitian case $Cent\,\k_{\Bbb C}$ is 1-dimensional and if $I$ is a set
indexing the $K_{\Bbb C}$-irreducible submodules $V^{i},\,i\in I$, of
$\p_{\Bbb C}$ then $I$ is a two element set and $$\p_{\Bbb C} =
\sum_{i\in I} V^{i}\eqno (3.36)$$ Also if $i\in I$ there exists
a linear isomorphism $\delta_{i}:Cent\,\k_{\Bbb C}\to \Bbb C$ such that
$$V^{i} = \{u\in \p_{\Bbb C}\mid ad\,x(u) = \delta_{i}(x)\,u\,\,\quad\forall x\in
Cent\,\k_{\Bbb C}\}\eqno (3.37)$$ In addition one has $$\delta_{i'} = -
\delta_{i}\eqno (3.38)$$ where $\{i,i'\} = I$. 

In case (1) we will say that $\g$ is of non-Hermitian type and in case (2) we will say
that $\g$ is of Hermitian type,

Recall $\g_{\psi} = Cent\,\n$ and  we have put $d = dim\,\g_{\psi}$. See \S 2.3 and
\S 2.4. \vs {\bf Proposition 3.11.} {\it If $d>1$ then $M$ is transitive on the
unit sphere $S$ (relative to ${\cal H}$) in $\g_{\psi}$. Furthermore $\g$ is of
Hermitian type if and only if (a) $d=1$ and (b)
$M$ operates trivially  on
$\g_{\psi}$.} \vs {\bf Proof.} The transitivity statement actually is stated in
Theorem 2.1.7 and the Remark which follows that theorem in [K-2]. However one need
not use that reference. We have already established what is needed to prove
transitivity on $S$ in the present paper. Indeed by (2.43) and (2.44)
one has that
$[\m,e]$ spans the
${\cal H}$ orthocomplement of $\Bbb R\,e$ in $\g_{\psi}$. This implies $Ad\,M (e)$ is
open in $S$. On the other it hand is also closed since $M$ is compact. Thus $Ad\,M (e)
= S$ if $d>1$ since $S$ is connected if $d>1$. 

In any case one has $$Ad\,G\, (e) = Ad\,K\,\,(\Bbb R^+ e)\eqno (3.39)$$ by the Iwasawa
decomposition $G = KAN$. But $Ad\,G\, (e)$ and hence also $Ad\,K\,(e)$ spans $\g$ by
the simplicity of $\g$. But if $\g$ is of Hermitian type then certainly
$Cent\,k\neq 0$ and if 
$0 \neq x\in Cent\,k$ the function $f(k) =(Ad\,k\,\,(e), x)$ on $K$ is constant,
real and non-zero. But then $-e\notin Ad\,M \,e$. Hence $d=1$ and $M$ operates
trivially on $\g_{\psi}$ by the transitivity statement. Conversely if $d=1$ and $M$
operates trivially on $\g_{\psi}$ then the line $\Bbb C\,e$ is stable under the
complex parabolic subalgebra $(\m + \a + \n)_{\Bbb C}$. But then $e$ generates a
complex irreducible $ad\,g_{\Bbb C}$-module, necessarily containing $\g$, and
hence equal to
$\g_{\Bbb C}$. Thus $\g_{\Bbb C}$ is simple and $e$ is the extremal weight vector
corresponding to $(\m + \a + \n)_{\Bbb C}$. But then $\g_{\Bbb C}$ is a spherical
$ad\,\g_{\Bbb C}$-module by the Cartan-Helgason theorem. Any non-zero spherical vector
must clearly lie in $\k_{\Bbb C}$ (e.g., by (3.37), (3.38) and the simplicity of
$g_{\Bbb C}$). But then $Cent\,k_{\Bbb C}\neq 0$. Thus $\g$ is of Hermitian type.
QED\vs {\bf Remark 3.12. } Note that as a consequence of Proposition 3.11 one has
the statement $$-e\notin Ad\,M(e)\,\,\,\iff \g\,\,\hbox{is of Hermitian type}\eqno
(3.40)$$ But one also has the statement $$-e\notin Ad\,G(e)\,\,\,\iff \g\,\,\hbox{is
of Hermitian type}\eqno (3.41)$$ Indeed if $\g$ is of non-Hermitian type then
$-e\in Ad\,M(e)$ by (3.40) and hence of course $-e\in Ad\,G(e)$. If $\g$ is of
Hermitian type let $0\neq x\in Cent\,k$. If $-e\in Ad\,G(e)$ then the function
$f$ in the proof of Proposition 3.11 must change signs by (3.39) and hence cannot be
constant. This proves (3.41).\vs Let $Proj(\g^*)$ be the real projective space defined
by $\g^*$ and let $P:\g^*-\{0\}\to Proj(\g^*)$ be the projectivization map. Let $C
=\gamma^{-1}(Cent\,\n)$ so that $\varepsilon\in C\s \g^*$ (see Proposition 2.5 and
(3.1)).\vs {\bf Theorem 3.13.} {\it There exists a unique closed $G$-orbit $\widehat
{Z}$ in
$Proj(\g^*)$ so that if $\widehat {Y}$ is any $G$-orbit in $Proj(\g^*)$ then
$\widehat {Z}$ is contained in the closure of $\widehat {Y}$. In particular
$$dim\,\widehat {Y}\geq dim\,\widehat {Z}\eqno (3.42)$$ and equality occurs in (3.41)
if and only if
$\widehat {Y} = \widehat {Z}$. 

Furthermore if $Z'\s \g^*$ is any non-zero $G$-coadjoint orbit then $Z'\s
P^{-1}(\widehat {Z})$ if and only if $Z'$ is of the form $Coad\, G\,(e')$ for $e'\in
C-\{0\}$. In fact, using the notation of Theorem 3.10 one has $$Z = P^{-1}(\widehat
{Z})\eqno (3.43)$$ in case $\g$ is of non-Hermitian type. In particular $Z =
-Z$ in the non-Hermitian case. If $\g$ is of Hermitian type let $J$ be an index set
parameterizing all
$G$-coadjoint orbits,
$Z^{j},\,j\in J$, in
$P^{-1}(\widehat {Z})$. Then $J$ is a 2-element set and one has $$\{Z,-Z\} =
\{Z^j\},\,\,j\in J\eqno (3.44)$$ so that $Z\neq -Z$ in the
Hermitian case. In either case one has $$dim\,Y\geq dim\,Z\eqno (3.45)$$ for any
non-zero coadjoint $G$-orbit $Y$.}\vs {\bf Proof.} Let $0\neq x\in \g$. Then the
span $\g_1$ of $Ad\,g(x)$ over all $g\in G$ is a non-zero ideal in $\g$ and hence
$\g_1= \g$. Thus there exists $g\in G$ such that the component, $x'$, of $Ad\,g(x)$
in
$\g_{\psi}$ with respect to the decomposition $$\g = \m + \a +\sum_{\varphi\in
\Delta}\g_{\varphi}$$ is not zero. But then, by Proposition 2.7, one has $$\lim_{t\to
+\infty} P(Ad\,(exp\,\,t\,x_{\psi})(x))= P(x')\eqno (3.46)$$ where we also let $P$
denote the projection map of $\g$ onto $Proj(\g)$. On the other hand if
$x''\in \g_{\psi}- \{0\}$ then $P(Ad\,G\,(x''))= P(Ad\,G\,(x'))$ by Proposition 3.11.
But $P(Ad\,G\,(x''))$ is compact since $P(x'')$ is fixed by the action of $AN$ and
hence $P(Ad\,G\,(x'')) = K\cdot P(x'')$. But $P^{-1}(P(Ad\,G\,(x''))$ is a single
$Ad\,G$-orbit by Proposition 3.11 and (3.40) in case $\g$ is of non-Hermitian type.
On the other hand it decomposes into a union of 2 distinct $Ad\,G$ orbits, one the
negative of the other, by Proposition 3.11 and (3.41), in case $\g$ is of Hermitian
type. Since $\gamma:\g^*\to \g$ is a $G$-equivalence the statements above
clearly carry over to $\g^*$. QED\vs {\bf Remark 3.14.} Because $\g$ is simple
and $Z$ is a $G$-coadjoint orbit there exists a Lie algebra injective homomorphism
$\g\to C^{\infty}(Z)$ (of course with respect to the Poisson algebra structure on
$C^{\infty}(Z)$). But then the symplectic isomorphism (3.30) implies that there exists
a Lie algebra injective homomorphism $$\g\to C^{\infty}(\widetilde {X})$$  By the
coadjoint orbit covering theorem and (3.45) note that $dim\,\widetilde {X}$ is
the smallest possible dimension of a symplectic manifold $W$ which admits an embedding
of
$\g$ as a Lie algebra of functions on $W$ under Poisson bracket. This points to an
interesting difference between
$(X,\omega)$ and the induced symplectic manifold $(\widetilde {X},\omega_{\widetilde
{X}})$ in that one cannot have a non-trivial Lie algebra homomorphism $$\g \to
C^{\infty}(X)$$ since $$dim\,X = dim\,Z\,-2\eqno (3.47)$$ by
(1.18).\vs If $V$ is a complex finite-dimensional
irreducible $K$ (and hence $K_{\Bbb C}$ module) let $X(V)\s \k^*$ be the
integral $K$-coadjoint orbit associated to $V$ by the Borel-Weil theorem. Thus if
$v\in V$ is an extremal weight vector there exists a Cartan subalgebra $\hh_{\Bbb C}$
of
$\k_{\Bbb C}$, which is the complexification of a Cartan subalgebra of $\k$, such that
$v$ is an  $\hh_{\Bbb C}$-weight vector for an extremal $\hh_{\Bbb C}$-weight
$\lambda$. We may regard $\lambda\in \k_{\Bbb C}^*$ where $\lambda\vert[\hh_{\Bbb
C},\,\k_{\Bbb C}]= 0$. Then $\nu = \lambda/2\,\pi\,i$ is real on $\k$ and $X(V)$ is
the coadjoint orbit of $\nu$. It is straightforward to show that $X(V)$ is
independent of the choice of $\hh_{\Bbb C}$. 

The embedding $\k\to \g$ induces, by transpose, a surjection $$\mu_{\g,\k}:\g^*\to
\k^*\eqno (3.48)$$ \vskip.5pc {\bf Remark 3.15.} Note that if $Y$ is a
$G$-coadjoint orbit then the restriction $\mu_{\g,\k}\vert Y$ is the moment map for
the action of $K$ on $(Y,\omega_Y)$ where $\omega_Y$ is the KKS-symplectic form on
$Y$. Indeed this follows immediately from (2.4) where $G$ replaces $K$ and $Y$
replaces
$X$. \vs {\bf Theorem 3.16.} {\it We use the notation of Theorem 3.10. One has
$$\mu_{\g,\k}(Z) = \Bbb R^{+}\,X\eqno (3.49)$$
 Moreover if 
$\g$ is of non-Hermitian type then $X=X(\p_{\Bbb C})$. Also $X=-X$ in
the non-Hermitian case. If $\g$ is of Hermitian type let $X^i= X(V^i)$ for
$i\in I$ (see (3.37)). Then $X \neq
-X$ and $$\{X,-X\} = \{X^i\},\,\,i\in I\eqno (3.50)$$ Furthermore (recalling
(3.44)) there exists a unique bijection $\alpha:I\to J$ such that
for $i\in I$, $$\mu_{\g,\k}(Z^{\alpha(i)}) =
\Bbb R^{+}\,X^{i}\eqno (3.51)$$ If $X = X^i$ then $Z = Z^{\alpha(i)}$. Finally,
adding to the symplectic isomorphism (3.30), for any $i\in I$, there exists a
symplectic $K$-diffeomorphism $$\beta_i:(\widetilde
{X^i},\omega_{\widetilde {X^i}})\to (Z^{\alpha(i)},\omega_{Z^{\alpha(i)}})\eqno
(3.52)$$ where $\omega_{Z^{\alpha(i)}}$ is the KKS form on the $G$-coadjoint orbit
$Z^{\alpha(i)}$.}\vs {\bf Proof.} Let $p:\g \to \k$ be the projection
corresponding to the Cartan decomposition $\g = \k + \p$. Obviously $$p\circ \gamma =
\gamma\circ \mu_{\g,\k}\eqno (3.53)$$ on $\g^*$. Recall (see \S 2.5) that $O$ is the
$G$-adjoint orbit of $e$. But $O = Ad\,K\,\,\Bbb R^+\,e$ by (2.56). Thus $f\in O$ if
and only if there exists $t\in \Bbb R^+$ and $k\in K$ such that $f =
t\,Ad\,k(e/\pi)$. But clearly
$p(e/\pi) = (e + \theta\,e)/(2\,\pi)$. That is $p(e/\pi) = h/(2\,\pi\,i)$ by (2.39).
Thus
$p(f) = t\,Ad\,k(h/(2\,\pi\,i))$. Applying $\gamma^{-1}$ to both sides one has
$$\mu_{\g,k}(\gamma^{-1}(f)) = t\,Coad\,\,k(\gamma^{-1}\,(h)/(2\,\pi\,i))$$ by
(3.53). But
$O = \gamma(Z)$ (see \S 3.1) so that if $\rho = \gamma^{-1}(f)$ then $\rho\in Z$ and
the most general element in $Z$ is of this form. Furthermore
$\gamma^{-1}\,(h)=
\lambda$ (see (2.49)) and
$\nu = \lambda/(2\,\pi\,i)$ (see \S 2.5). Thus $\mu_{\g,k}(\rho) =
t\,Coad\,k(\nu))$. Since $k\in K$ and $t\in \Bbb R^+$ are arbitrary and $X$ is
the $K$-coadjoint orbit of $\nu$ this proves (3.49). 

Now since $B\vert \p_{\Bbb C}$ is non-singular it follows that $\p_{\Bbb C}$ as a
$K$-module is self-contragedient. In particular $$\hbox{$\mu$ is a weight of $\p_{\Bbb
C}\,\,\iff\,\,-\mu$ is a weight of $\p_{\Bbb C}$}\eqno (3.54)$$  

Recall the choice of $X$ in \S 2.5. By definition $$ X= X(V_{\lambda})\eqno 
(3.55)$$ where $V_{\lambda}$ is the $K$-irreducible submodule of $\p_{\Bbb C}$
defined in \S 2.5 with extremal weight vector $v$ (see (2.48)). Assume $\g$ is
of non-Hermitian type. Then $\p_{\Bbb C}$ is $K$-irreducible so that
$V_{\lambda} =
\p_{\Bbb C}$. This proves $X = X(\p_{\Bbb C})$. But then $X = -X$ by (3.54). 

Now assume that $\g$ is of Hermitian type. Then $V_{\lambda} = V^i$ for some $i\in
I$. But $\mu\vert Cent\,k = \delta_i\vert Cent\,\k$ for any weight $\mu$ of
$V_{\lambda}$ by (3.37). Thus $X \neq -X$ and one easily has (3.50) by (3.54). But
(3.49) implies
$\mu_{\g,\k}(-X) = \Bbb R^+(-Z)$. But then (3.51) follows from (3.44) and (3.49). 

The choice of $e\in \g_{\psi}$ was arbitrary (see \S 2.3) subject only to the
condition that $\{e,e\}= 1$. We can therefore replace $e$ by $-e$. In that case $h$
would be replaced by $-h$ and $\lambda$ by $-\lambda$. Consequently $X$ is replaced
by $-X$ and $Z$ by $-Z$. Thus the argument which leads to the symplectic
$K$-isomorphism (3.30) yields a symplectic $K$-diffeomorphism $(\widetilde
{-X},\omega_{\widetilde {-X}})\to (-Z,\omega_{-Z})$. This implies (3.52). QED\vs
3.4. Even though $\g$ is simple the complexification $\g_{\Bbb C}$ may not be
simple. Indeed this is the case if $\g$ itself were complex but its complex
structure is ignored. We will say that $\g$ is 
$O_{min}$-split if $\g_{\Bbb C}$ is a simple complex Lie algebra and
$e\in O_{min}$ where $O_{\min}$ is the minimal nilpotent orbit in
$\g_{\Bbb C}$. Note that by the transitivity statement in Proposition 3.11 and
Theorem 3.13 the definition of $O_{min}$-split depends only on $\g$ and is, in 
particular, independent of the choice of
$e$. If $\g$ is $O_{min}$-split then $$dim\,Z = 2h^{\vee} - 2\eqno (3.56)$$ where
$h^{\vee}$ is the dual Coxeter number of $\g_{\Bbb C}$. This is clear since one knows
$dim_{\Bbb C}O_{\min} = 2h^{\vee} - 2$. But by the irreducibility of the adjoint
representation of $\g_{\Bbb C}$ one has $O_{\min} = Ad\,\g_{\Bbb C}(e)$. However $O =
Ad\,\g(e)$ (see
\S 2.5). But then (3.56) follows by computing the dimension of the respective
tangent spaces at $e$. See \S 3.1 for the definition of $Z$. 
\vs {\bf Theorem 3.17.}  {\it The simple Lie algebra $\g$ is $O_{min}$-split if and
only if $$dim\,Cent\,\n =1\eqno (3.57)$$}\vs {\bf Proof.} If $dim\,Cent\,\n=1$ then
$\Bbb C\,e$ is stabilized, under the adjoint represention, by the complex parabolic
subalgebra $(\m + \a +\n)_{\Bbb C}$. But then $e$ generates an irreducible $\g_{\Bbb
C}$-module under the adjoint representation of $\g_{\Bbb C}$. But this module
contains $\g$ by the simplicity of $\g$. Hence the module equals $\g_{\Bbb C}$ so
that $\g_{\Bbb C}$ is simple. Furthermore $e\in O_{\min}$ since $e$ is an extremal
weight vector of this module. Thus $\g$ is $O_{min}$-split. Conversely assume that
$\g$ is $O_{min}$-split. Then $\u_{\psi}$ (see \S 2.3) is conjugate to the TDS of a
highest root vector in $\g_{\Bbb C}$ by Corollary 3.6 in [K-3]. Since the root in
question is long the multiplicity of the eigenvalue 2 of $ad\,x_{\psi}$ is 1. Thus
$d=1$ by Proposition 2.7. That is, $dim\,Cent\,\n = 1$. See \S 2.4. QED\vs 
{\bf Examples. 3.18.} By Theorem 3.17 $\g$ is $O_{min}$-split if (a) $\g$ is
split (so that all restricted root spaces are 1-dimensional), (b)
$\g_{\Bbb C}$ is simple and
$\g$ is a quasi-split form of $\g_{\Bbb C}$ (since $\n_{\Bbb C}$ is a maximal
nilpotent Lie algebra of
$\g_{\Bbb C}$), (c) 
$\g$ is of hermitian type (by Proposition 3.11). \vs Note that if $\g$ is
$O_{min}$-split then $$dim\,X = 2h^{\vee} - 4\eqno (3.58)$$ by (3.47) and (3.56). In
general
$X
\cong K/K_{\nu}$ can be given the structure of a partial complex flag manifold (see
(2.16)). If
$\g$ is $O_{min}$-split much more can be said. \vs {\bf Proposition 3.19.} {\it
If $\g$ is $O_{min}$-split (e.g., if $\g$ is split) then $X$ is not only a
$K$-symmetric space but in fact $X$ is a Hermitian symmetric space.}\vs {\bf Proof.}
Recall the notation of \S 2.4. If $\g$ is $O_{min}$-split then by Theorem 3.17
the number $d_{\k}$ in Theorem 2.11 is 0. Thus, since $x_{\psi}$ and $h$ are
conjugate, the eigenvalues of $ad\,h\vert \k_{\Bbb C}$ are in the set $\{1,0,-1\}$.
Since
$(\k_{\nu})_{\Bbb C}$ is the $ad\,h\vert \k_{\Bbb C}$-eigenspace for the eigenvalue
$0$ the result is immediate. QED\vs {\bf Remark 3.20}. Assume that $\g$ is
$O_{min}$-split and $\g$ is of non-Hermitian type so that $\k_{\Bbb C}$ is
semisimple. Let $K'$ be a non-compact real form of $K_{\Bbb C}$ having $K_{\nu}$ as a
maximal compact subgroup so that
$K'/K_{\nu} = X'$ is the non-compact symmetric dual to the compact symmetric space
$X$. But now $X=-X$ by Theorem 3.16 so that $h$ and $-h$ are $K$-conjugate.
One readily shows that this implies that $h$ lies in a TDS of $\k_{\Bbb C}$ so that
not only is $X'$ a complex bounded domain but in fact $X'$ is a tube
domain. See [K-W]. In particular by the Kantor-Koecher-Tits theory
$X'$ corresponds to a formally real Jordan algebra $J(X)$. See \S 5 in [K-S]
for a classification of the simple formally real Jordan algebras and the corresponding
tube domains.  One then has
$$dim\,X = 2\,dim\,J(X)\eqno (3.59)$$

If $\g$ is a split form of any one of 5 exceptional simple Lie
algebras, $\g$ is non-Hermitian so that Remark 3.20 applies to $\g$. But a minimal
(dimensional) symplectic realization of $\g$ as functions on a symplectic manifold is
achieved when the manifold is the induced symplectic $\widetilde {X}$ of a coadjoint
orbit of a compact Lie group $K$. See Remark 3.14. The group $K$ turns out to be
classical in all 5 cases so that the exceptional Lie algebras $\g$ emerge
symplectically from the symplectic induction of a classical coadjoint orbit. The
table below contains the relevant information. The cases of
$E_6,\,E_7$ and
$E_8$ are taken from [B-K]. To avoid complicated notation the compact groups
listed below are correct only up to finite coverings. The symbol $H(n)/F$
denotes the Jordan algebra of $n\times n$ Hermitian matrices over the field $F$. The
compact form of $Sp(2n,\Bbb C)$ will be denoted simply by $Sp(2n)$.

$$\matrix{\underline {\g_{\Bbb C} \,type}&\underline {K}&\underline
X &\underline {dim\,X}&\underline {J(X)}\cr 
G_2 & SU(2)\times SU(2)& P_1(\Bbb C)\times P_1(\Bbb C)&4& \Bbb R \oplus \Bbb R\cr F_4
&SU(2)\times Sp(6) & P_1(\Bbb C)\times (Sp(6)/U(3)) &14&\Bbb R \oplus (H(3)/{\Bbb
R})\cr E_6&Sp(8)& Sp(8)/U(4) &20& H(4)/{\Bbb R}\cr E_7&SU(8)&SU(8)/(SU(4)\times
SU(4)\times U(1))& 32 & H(4)/{\Bbb C}\cr E_8& Spin(16)&Spin(16)/U(8)& 56 &
H(4)/{\Bbb H}\cr}$$\vskip .5pc

\centerline{\bf References}\vskip 1.3pc
\parindent=42pt
\rm
\item {[B-K]} R. Brylinski and B. Kostant, Lagrangian Models of
Minimal Representations of $E_6$, $E_7$, and $E_8$, {\it Progress
in Math. Birkhauser Boston}, Vol. {\bf 131}(1995), 13-53
\item {[K-1]} B. Kostant, Quantization and Unitary representations,
{\it Lecture Notes in Math. Springer-Verlag}, Vol. {\bf 170}(1970), 87-207
\item {[K-2]} B. Kostant, On the existence and irreducibility of
certain series of representations, {\it Lie groups and Their
Representations}, edited by I. M. Gelfand, John Wiley \& Sons,
1971, 231-329
\item {[K-3]} B. Kostant, The principal three-dimensional subgroup
and the Betti numbers of a complex simple Lie group, {\it Amer. J.
Math.} {\bf 81}(1959), No. 4, 973-1032
\item {[K-S]} B. Kostant and S. Sahi, The Capelli Identity, Tube
Domains, and the Generalized Laplace Transform, {\it Adv. in
Math.}, {\bf 87}(1991), 71-92
\item {[K-W]} A. Koranyi and J. Wolf, Realization of Hermitian
Symmetric spaces as generalized half-planes. {\it Ann. of Math.}
{\bf 81}(1965), 265-288 
\item {[S]} J. Sekiguchi, Remarks on real nilpotent orbits of a
symmetric pair, {\it J. Math. Soc. Japan}, {\bf 39}(1987), 127-138
\item {[V]} M. Vergne, Instantons et correspondance de
Kostant-Sekiguchi, {\it C. R. Acad. Sci. Paris}, {\bf
320}(1995), S\'erie 1, 901-906

\smallskip
\parindent=30pt
\baselineskip=14pt
\vskip 1.9pc
\vbox to 60pt{\hbox{Bertram Kostant}
      \hbox{Dept. of Math.}
      \hbox{MIT}
      \hbox{Cambridge, MA 02139}}\vskip 1pc

      \noindent E-mail kostant@math.mit.edu

\end